\documentclass{cmslatex}
\usepackage[paperwidth=7in, paperheight=10in, margin=.875in]{geometry}

\usepackage{amstext}
\usepackage[numbers,sort&compress]{natbib}
\usepackage{graphicx}
\usepackage{amsmath, caption, subcaption, amssymb}
\usepackage{color}
\usepackage{todonotes}
\usepackage{float}
\usepackage{hyperref}
\hypersetup{
    colorlinks,
	linkcolor={red!50!black},
    citecolor={blue!70!black},
	citecolor={blue},
    urlcolor={blue!80!black}
}

\allowdisplaybreaks

\newcommand{\figref}[1]{Figure~\ref{#1}}
\newcommand{\secref}[1]{Section~\ref{#1}}
\newcommand{\appref}[1]{Appendix~\ref{#1}}

\renewcommand{\theequation}{\arabic{section}.\arabic{equation}}

\begin{document}
\title{Coalescing particle systems and applications to nonlinear Fokker-Planck equations\thanks{}}


\author{Gleb Zhelezov\thanks{Department of Mathematics, University of Arizona, 617 N Santa Rita Ave, PO Box 210089, Tucson, AZ 85718,  (gzhelezov@math.arizona.edu) \url{http://math.arizona.edu/\~gzhelezov}} \and{Ibrahim Fatkullin\thanks{Department of Mathematics, University of Arizona, 617 N Santa Rita Ave, PO Box 210089, Tucson, AZ 85718  (ibrahim@math.arizona.edu) \url{http://math.arizona.edu/\~ibrahim}}}}




\pagestyle{myheadings} \markboth{COALESCING PARTICLE SYSTEMS}{G. ZHELEZOV AND I. FATKULLIN} \maketitle

\begin{abstract}
We study a stochastic particle system with a logarithmically-singular inter-particle interaction potential which allows for inelastic particle collisions. We relate the squared Bessel process to the evolution of localized clusters of particles, and develop a numerical method capable of detecting collisions of many point particles without the use of pairwise computations, or very refined adaptive timestepping. We show that when the system is in an appropriate parameter regime, the hydrodynamic limit of the empirical mass density of the system is a solution to a nonlinear Fokker-Planck equation, such as the Patlak-Keller-Segel (PKS) model, or its multispecies variant. We then show that the presented numerical method is well-suited for the simulation of the formation of finite-time singularities in the PKS, as well as PKS pre- and post-blow-up dynamics. Additionally, we present numerical evidence that blow-up with an increasing total second moment in the two species Keller-Segel system occurs with a linearly increasing second moment in one component, and a linearly decreasing second moment in the other component.
\end{abstract}



\section{Introduction}
\subsection{Background}

The connection between systems of interacting particles and kinetic-type PDEs was first investigated by Kac in his study of the motion of a tagged molecule in a bath of identical molecules \cite{kac1959probability}, which arose as a simplified model of a Maxwellian gas \cite{mckean1966speed}. This work introduced the property of ``propagation of chaos'': as the number of molecules tends to infinity, the $N$-particle probability densities are well-approximated by the product of single particle marginals.

The connection between such processes and nonlinear parabolic equations, such as Boltzmann's equation or Burgers' equation, was then elaborated by McKean \cite{mckean1966class}. This line of research has continued since, and much more is now known about the duality between these processes and parabolic PDEs \cite{sznitman1991topics}. In particular, particle-based numerical methods have been developed for the solution of such PDEs \cite{carrillo2014derivation} using the methods of ``mean field Monte Carlo.'' The solutions to these PDEs are approximated by the empirical density of $N$-particle systems. As the number of particles tends to infinity, such approximations become exact by the propagation of chaos property.

Rigorously proving propagation of chaos for particle systems with singular interaction coefficients is challenging, and has only been carried out in a few special cases, e.g. \cite{jabin2011particles}. One PDE associated with a logarithmically-singular particle system is the Patlak-Keller-Segel chemotaxis model (PKS) \cite{keller1970initiation, patlak1953random}, which is reviewed extensively in \cite{horstmann20031970, horstmann20041970}. Despite the lack of a propagation of chaos result, the PDE has been numerically approximated using the associated particle system in several works, initially in \cite{haskovec2009stochastic, haskovec2011convergence} and later in \cite{fatkullin2013study}. Various properties of the PKS, such as the formation of Dirac singularities in finite time \cite{blanchet2006two}, as well as interaction of singularities post-blow-up \cite{velazquez2004pointsI, velazquez2004pointsII, dolbeault2009two}, can either be shown to be true in the particle system, or have considerable numerical evidence for their existence. Recent advacements in understanding this particle system include partial existence and uniqueness results for solutions to the subcritical (small mass) particle system \cite{fournier2015stochastic, cattiaux20162}, and convergence of the empirical density of a similar particle system to the solution of a modified PKS system \cite{budhiraja2016uniform}.

Singular interaction coefficients in the PKS particle system  allow for particle collisions, and some type of regularization must be introduced in order to propagate the particle system past the first collision time. In \cite{haskovec2009stochastic}, semi-deterministic heavy particles absorb light particles. In \cite{fatkullin2013study}, collided particles are forced to move in unison due to a mean field. Broadly speaking, the two works take two different approaches to simulating the regularizations of the PKS derived in \cite{dolbeault2009two} and \cite{velazquez2004pointsI, velazquez2004pointsII}. The first work simulates the singular limit of the system, whereas the second work simulates the system with an effectively regularized Green's function.

In \cite{haskovec2009stochastic}, heavy particles corresponding to singularities in the PDE must be prescribed a priori and cannot arise as the result of a collision of many light particles. On the other hand, particles do not truly collide in \cite{fatkullin2013study}, and the deterministic system approximated is closer to the one given in \cite{velazquez2004pointsI, velazquez2004pointsII}, where singularities are replaced with regions of high density. In this work, we develop criteria for particle coalescence of particles of arbitrary masses, based on analytical estimates of exit times of the squared Bessel process. In this context, the particle system in \cite{haskovec2009stochastic} can be viewed as the limit of the particle system in \cite{fatkullin2013study} with collisions, as the number of particles tends to infinity.

\subsection{Outline} We introduce a coalescing particle system with nonuniform particle masses and a logarithmic interaction kernel. Using estimates on the system's second moment, we derive a criterion for a finite-time collision of the entire particle system. We then motivate the mass-dependence of the diffusion coefficient of a particle, and approximate the time evolution of a localized subsystem's second moment. We then show that the hydrodynamic limit of such a system is the multispecies Patlak-Keller-Segel system, of which the PKS is a special case. Finally, we present a numerical method implementing many-body collisions and coalescence events, which is generally applicable to PDEs of the form
\begin{align}
\label{eq:nonlinear_FP}
\begin{cases}
	\partial_t \rho_1 &= \nabla \cdot (\mu_1 \nabla \rho_1 - \chi \rho_1 \nabla c), \\
	&\vdots \\
	\partial_t \rho_K &= \nabla \cdot (\mu_K \nabla \rho_K - \chi \rho_K \nabla c), \\
	\mathcal{L}c &= -(\rho_1 + \cdots + \rho_K),
\end{cases}
\end{align}
where
\begin{equation}
	\label{eq:elliptic_operator}
	\mathcal{L}c(x,t) = \nabla \cdot \left( G(x) \nabla c(x,t) \right) + F(x,c)
\end{equation}
is an elliptic operator with a fundamental solution $V$ which has a logarithmic singularity. As an application, we apply it to the planar case with decaying (radiative) boundary conditions and $\mathcal{L} = \Delta$, though the method is equally applicable to bounded domains with Neumann boundary conditions. This special case is the planar PKS system, some properties of which we describe in \secref{sec:PKS_intro}, and whose measure-valued solutions we describe in \secref{sec:PKS}. We also apply the numerical method to investigate blow-up in the components of the multispecies PKS.

\subsection{The coalescing particle system}
We study the $N$-particle systems described by the following equations
\begin{equation}
	\label{eq:particle_dynamics}
	dX^{(n)}_t = - \chi \frac{\partial}{\partial X^{(n)}_t} \sum_{\substack{i = 1 \\ i \neq n}}^N m_i V(X^{(n)}_t, X^{(i)}_t) dt + \sqrt{\frac{2\tilde{\mu}}{m_n}} dW^{(n)}_t.
\end{equation}
Each particle has some mass $m_n$ and position $X^{(n)}_t \in \mathbb{R}^2$. The total mass is $M = \sum_i m_i $, and $\chi, \tilde{\mu} > 0$ are parameters. The processes $W^{(n)}_t$ are independent Wiener processes.

The particle system in \eqref{eq:particle_dynamics} is related to the PDE in \eqref{eq:nonlinear_FP} when $V$ is the fundamental solution of $\mathcal{L}$, e.g. if $\mathcal{L} = \Delta$ or $\mathcal{L} = \Delta - k^2$, we have
\begin{subequations}
\begin{align}
	\label{eq:PKS_greens_function_1}
	V(x,y) &= \frac{1}{2\pi} \ln |x-y|, \\
	\label{eq:PKS_greens_function_2}
	V(x,y) &= -\frac{1}{2\pi} K_0(k|x-y|),
\end{align}
\end{subequations}
where $K_0$ is the modified Bessel function of the second kind. When $m_n = M/N$ and $\tilde{\mu} = \mu M / N$, the empirical mass density of the particle system with \eqref{eq:PKS_greens_function_1} approximates the PKS, and the particle system with \eqref{eq:PKS_greens_function_2} is the one given in \cite{fatkullin2013study}.

The dynamics prescribed in equation \eqref{eq:particle_dynamics} allows for particle collisions provided that $V$ has logarithmic or stronger singularities. In this case, the SDE must be augmented with proper boundary conditions prescribing behavior when at least two particles' coordinates are identical. Well-posedness and uniqueness results for these types of SDEs have not been rigorously established. We proceed formally, considering inelastic collisions: colliding particles merge into a single particle which absorbs their total mass.

\subsection{\label{sec:PKS_intro}Properties of the Patlak-Keller-Segel system} Since many of the applications of this work are related to the PKS, we give a short overview of its definition and properties here.

The PKS is prescribed by the following system of PDEs:
\begin{align}
\begin{cases}
	\partial_t \rho &= \nabla \cdot ( \mu \nabla \rho - \chi \rho \nabla c ), \\
	\Delta c &= -\rho,
\end{cases}
\end{align}
and models a biological system consisting of amoeba, which spread across the plane with mass density $\rho(x,t)$ and produce a chemical (``chemoattractant'') of concentration $c(x,t)$. On average, amoeba diffuse in space with diffusivity $\mu$ and drift in the direction of $\nabla c$ with speed $\chi | \nabla c|$. The chemoattractant diffuses instantly. The boundary condition $\rho(x,t) \rightarrow 0$ as $|x| \rightarrow \infty$ is enforced, and mass is conserved: $\int \rho(x,t)dx = M$.

This system has been investigated extensively in the literature \cite{horstmann20031970, horstmann20041970}, often in connection with the property that when
\begin{equation}
	M > 8\pi\mu / \chi,
\end{equation}
solutions form singularities in finite time, and when
\begin{equation}
	M < 8\pi\mu / \chi,
\end{equation}
solutions are global in time \cite{blanchet2006two}. In the former case, an upper bound for the singularity formation time $T$ may be given as
\begin{equation}
	T < \frac{2\pi F(0)}{(\chi M - 8\pi\mu)M},
\end{equation}
where
\begin{equation}
	F(t) = \int_{\mathbb{R}^2}|x|^2\rho(x,t)dx
\end{equation}
and $F(0)$ is the system's initial second moment \cite{blanchet2006two}.

\section{Collisions and post-collision dynamics}
\subsection{Overview}
Let us first carry out a moment-based computation for finding a criterion which predicts whether a particle system will coalesce into a single particle in finite time. Similar to the PKS mass criterion, this criterion only depends on the total mass of the system and the number of particles, and is otherwise independent of the distribution of particles in the plane. We then motivate the mass dependence of the diffusion coefficient of the newly created particle. Finally, we derive an approximate equation for the dynamics of the second moment of an isolated cluster of particles.

\subsection{\label{sec:collision_criterion}Collision criterion for the full system}
\begin{figure}
	\captionsetup{width=\linewidth}
	\centering
	\includegraphics[scale=0.8]{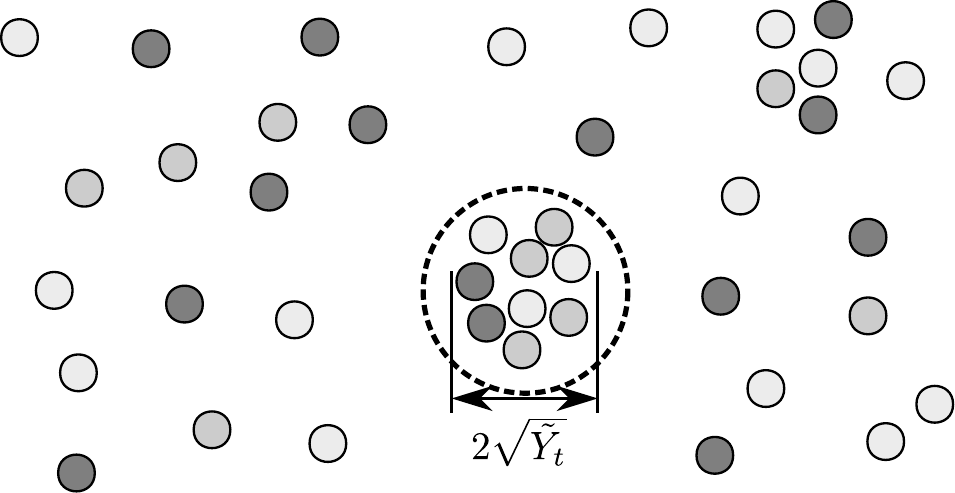}	\caption{An $N$-particle system with a tightly-clustered $N'$-particle subsystem. The particles inside the dashed circle correspond to particles with indices $1, \ldots, N'$, and the rest of the particles correspond to $N'+1, \ldots, N$. Several colors are used to emphasize that the point particles are of different masses.}
	\label{fig:formed_cluster}
\end{figure}
Consider an $N$-particle system with masses and $V$ given as in \eqref{eq:PKS_greens_function_1}. The dynamics of the $n$th particle are then prescribed by
\begin{equation}
	dX^{(n)}_t = -\frac{\chi}{2\pi} \sum_{i \neq n} m_i \frac{X^{(n)}_t - X^{(i)}_t}{\left| X^{(n)}_t - X^{(i)}_t \right|^2} dt + \sqrt{\frac{2\tilde{\mu}}{m_n}} dW^{(n)}_t.
\end{equation}
To quantify the size of the system, consider its second moment
\begin{equation}
	Y_t = \frac{1}{2M^2} \sum_{i,j} m_i m_j \left| X^{(i)}_t - X^{(j)}_t \right|^2.
\end{equation}
By the positivity of $Y_t$, showing the total collision of the particles in finite time is equivalent to showing that $Y_T = 0$ for some $T < \infty$.

It can be shown (by an application of Ito's lemma) that
\begin{equation}
	\label{eq:dY_t}
	dY_t = \alpha dt + 2 \beta \sqrt{Y_t} dW_t
\end{equation}
where
\begin{align}
	\label{eq:alpha_beta}
	\alpha &=  \frac{4\tilde{\mu}(N-1)}{M} - \frac{\chi M}{2\pi} \left(1 - \sum_j \left( \frac{m_j}{M} \right)^2 \right), \quad \beta = \sqrt{\frac{2\tilde{\mu}}{M}},
\end{align}
and
\begin{equation}
	\label{eq:Second_moment_driving_process}
	dW_t = \frac{1}{(M)^{3/2} \sqrt{Y_t}} \sum_{i,j=1}^{N} m_j \sqrt{m_i} \left( X^{(i)}_t - X^{(j)}_t \right) \cdot dW^{(i)}_t
\end{equation}
is a Wiener process by the L\'{e}vy characterization. We stress that expression \eqref{eq:dY_t} is only valid between collision events, as $\alpha$ depends on the total number of particles and their masses, and must therefore be updated after each collision. Rescaling time as $t \rightarrow t/\beta^2$ and setting $\tilde{Y}_t = Y_{\beta^2 t}$, we get
\begin{equation}
	\label{eq:Bessel_process}
	d\tilde{Y}_t = 2(\nu + 1) dt + 2 \sqrt{\tilde{Y}_t} dW_t,
\end{equation}
where $\nu = \frac{\alpha}{2\beta^2} - 1$. In terms of our original constants, $\nu$ is given by
\begin{equation}
	\label{eq:Bessel_index}
	\nu(m_1, m_2, \cdots, m_N) = (N - 2) - \frac{\chi M^2}{8 \pi \tilde{\mu}} \left( 1 - \sum_{j} \left(\frac{m_j}{M}\right)^2 \right).
\end{equation}
Equation \eqref{eq:Bessel_process} describes a squared Bessel process with index $\nu$. Its boundary behavior at $\tilde{Y} = 0$ depends on its index \cite{revuz1999continuous, karlin1981second}:
\begin{enumerate}
	\item When $\nu \in [0, +\infty)$, the origin is an entrance boundary, and $\tilde{Y}_t > 0$ a.s. for all $t > 0$
	\item When $\nu \in (-1, 0)$, the origin is a regular boundary, and the behavior of the process at this point must be defined (e.g. absorbing boundary, reflective boundary)
	\item When $\nu \in (-\infty, -1]$, the origin is an absorbing boundary which is hit in finite time
\end{enumerate}
It then follows that a full, simultaneous collision of all the particles may occur if
\begin{equation}
	\label{eq:negative_index}
	\nu(m_1, \ldots, m_N) < 0.
\end{equation}
When $\nu \in (-1,0)$, we may choose the collision, which we call ``soft,'' to be fully inelastic, or fully elastic. Similarly, when $\nu~\in~(-\infty, -1]$, only an inelastic collision may occur.

The above is not sufficient for describing all collisions in the system. For instance, we expect the associated singular forces to force the subsystem inside the dashed line in \figref{fig:formed_cluster} to inelastically collide earlier than the full system. We will approximate the evolution of the second moment of such a colliding subsystem in \secref{sec:subsystem_second_moment}, but already note here that a localized colliding subsystem's second moment may be approximated as a separate squared Bessel process that's independent of the particles not participating in the collision. As shown in \appref{app:Subtraction_formula}, the indices of the squared Bessel processes corresponding of the full system pre- and post-collision, and the index of the colliding subsystem, are related via a subtraction formula: if $\nu_i$ is the index of the full system described in \figref{fig:formed_cluster}, $\nu_f$ is the index of the same system after the particles inside the dashed lines coalesce, and $\overline{\nu}$ is the index of the subsystem inside the dashed line, then
\begin{equation}
	\label{eq:index_difference}
	\nu_f - \nu_i = -\left( \overline{\nu} + 1 \right).
\end{equation}
From \eqref{eq:index_difference} we see that hard collisions, except in the critical $\overline{\nu}=-1$ case, always increase the system's overall index, and soft collisions increase the system's overall index.

To see the effect of this index change on the full system, let $\tau$ be the first hitting time of the origin for the SDE given in \eqref{eq:dY_t}. This hitting time has the inverse gamma distribution \cite{makarov2010exact}, 
\begin{equation}
	\label{eq:FHT_rv}
	\tau \asymp \tilde{\mu} \frac{Y_0}{U},
\end{equation}
where $U \sim G(|\nu|, 1)$ is distributed according to the gamma distribution with shape parameter $|\nu|$ and rate parameter $1$. 

Intuitively, we see that increasing the index implies that a system contracts at a slower rate, and that a system with only hard inelastic collisions contracts at a slower rate after each collision (e.g. as in \figref{fig:PKS_simulation_zoom}). Furthermore, we expect many systems which can experience soft inelastic collisions to behave similarly, as a localized subsystem with an index $\overline{\nu} \in (-1,0)$ has a low probability of undergoing a collision in a time step (e.g. $\tau$ only has an expected value when $\nu < -1$), and may attract a sufficient number of additional particles into its aggregate to force the aggregate to experience a hard collision instead. Since in this work we will primarily focus on the large particle case, we prescribe that all collisions---soft and hard (i.e. $\nu < 0$)---are inelastic.

We remark that the formula for the time derivative of the second moment of the PKS also only gives an upper bound for the formation of a singularity, since for a system of total mass $M$ greater than the system's critical mass $M_c$, a second moment equal to zero implies the formation of a singularity of total mass $M > M_c$. However, singularities in the radially-symmetric PKS form with a mass of exactly $M_c$ \cite{herrero1996singularity, velazquez2002stability}, after which the time derivative of the second moment changes \cite{velazquez2004pointsI}.

\subsection{Post-collision dynamics}
\begin{figure}
	\captionsetup{width=\linewidth}
	\centering
	\begin{subfigure}[t]{0.4\textwidth}
		\frame{\includegraphics[width=\textwidth]{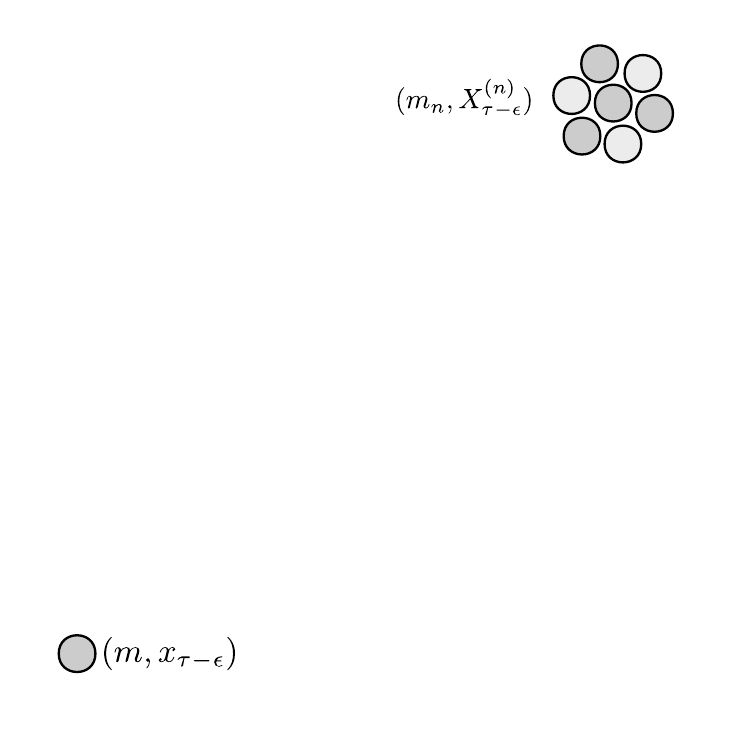}}
		\caption{An aggregate of particles, a~moment before coalescence.}
		\label{fig:equivalent_precoalesced}
	\end{subfigure}
	\begin{subfigure}[t]{0.4\textwidth}
		\frame{\includegraphics[width=\textwidth]{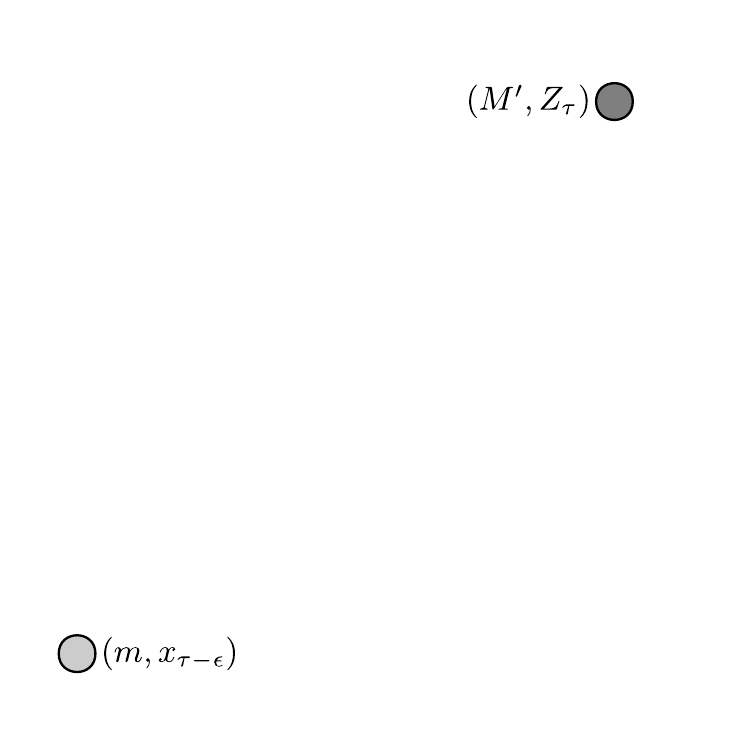}}
		\caption{Aggregate coalesces into one particle of mass $M'$.}
		\label{fig:equivalent_coalesced}
	\end{subfigure}
	\caption{As $\epsilon \rightarrow 0$, the bottom-left particle should experience the same drift in \eqref{fig:equivalent_precoalesced} and \eqref{fig:equivalent_coalesced}.}
	\label{fig:equivalent_particle}
\end{figure}
The dynamics of the coalescing diffusion system, given by \eqref{eq:particle_dynamics}, are undefined at times when there exist two indices $i$ and $j$ such that $X^{(i)}_t = X^{(j)}_t$. If we prescribe that collisions only occur inelastically, we can propagate the system past collision times by coarsening the system: that is, by replacing each collided aggregate of particles with a single particle of the same mass as the aggregate. Let us now show the diffusion coefficient of the newly-created particle is inversely-dependent on the square root of the mass, as given in \eqref{eq:particle_dynamics}.

Consider an $N' + 1$ particle system, with the first $N'$ particles positioned in a tight, pre-coalesced cluster at $X^{(n)}_t$ with masses $m_n$ totalling to $M'$, and the last particle located far away at $x = X^{(N'+1)}_t$ with mass $m = m_{N'+1}$, as in \figref{fig:equivalent_precoalesced}. In general, the diffusion coefficient of a particle may be given as a function of the particles mass, $\sigma_n = \sigma(m_n)$. Let $\tau$ denote the time at which the first $N'$ particles coalesce at $Z_\tau$, and fix $0 < \epsilon \ll \tau$. Then
\begin{equation}
	\label{eq:x_precoalescence_dynamics}
	dx_{\tau - \epsilon} = -\chi \sum_{i=1}^{N'} m_i \frac{\partial V}{\partial x} \left(x_{\tau-\epsilon},X^{(i)}_{\tau - \epsilon} \right)  dt + \sigma(m_{N'+1}) dW^{(N'+1)}_{\tau - \epsilon}.
\end{equation}
At the moment the first $N'$ particles coalesce, the system becomes a two-particle system, and so
\begin{equation}
	\label{eq:x_postcoalescence_dynamics}
	dx_\tau = -\chi M' \frac{\partial V}{\partial x} \left(x_{\tau},Z_{\tau} \right) dt + \sigma(m_{N'+1}) dW^{(N'+1)}_\tau.
\end{equation}
Let us assume the particle at $x_t$ should not experience an abrupt discontinuity in its drift at the moment of coalescence, i.e. we want $dx_{\tau - \epsilon} \rightarrow dx_\tau$ as $\epsilon \rightarrow 0^+$. Equating the right hand sides of \eqref{eq:x_precoalescence_dynamics} and \eqref{eq:x_postcoalescence_dynamics} as $\epsilon \rightarrow 0^+$ and using the property that $X^{(n)}_{\tau - \epsilon} \rightarrow Z_\tau$ for all $n \leq N'$, we get
\begin{equation}
	\label{eq:COM_limit}
	Z_\tau = \lim_{\epsilon \rightarrow 0^+} \frac{1}{M'} \sum_{i=1}^{N'} m_i X^{(i)}_{\tau - \epsilon},
\end{equation}
meaning the $N'$ particles must coalesce at the center of mass of the subsystem. This suggests that the diffusion coefficient of the newly-created particle positioned at $Z_\tau$ should be the same as the diffusion coefficient of the center of mass process of the first $N'$ particles for $t < \tau$. By the independence of the processes $W^{(i)}_t$ for $1 \leq i \leq N'$ and the definition of the center of mass inside the limit on the right hand side of \eqref{eq:COM_limit}, we get
\begin{equation}
	\sigma(M') = \frac{1}{M'} \sqrt{\sum_{i=1}^{N'} m_i^2 (\sigma(m_i))^2},
\end{equation}
or equivalently,
\begin{equation}
	(M')^2 \left(\sigma(M'\right))^2 = \sum_{i=1}^{N'} m_i^2 (\sigma(m_i))^2.
\end{equation}
Since $M' = \sum m_i$, it follows that $f(x) = x^2 (\sigma(x))^2$ must be additive, i.e. satisfies Cauchy's functional equation,
\begin{equation}
	f(x) = f(x) + f(y).
\end{equation}
Under the physically relevant assumption that $f$ is continuous, solutions to this functional equation must be linear \cite{kuczma2009introduction}. We therefore get
\begin{equation}
	\sigma(m) = \sqrt{\frac{2\tilde{\mu}}{m}},
\end{equation}
as in the dynamics given in the beginning of the work in \eqref{eq:particle_dynamics}.

By the same reasoning, we expect $Z_t$ to be driven by the weighted noise of the center of mass, $W^{(cm)}_t$, given by
\begin{equation}
	\label{eq:COM_noise}
	W^{(cm)}_t = \frac{1}{\sqrt{M'}} \sum_{i=1}^{N'} \sqrt{m_i} W^{(i)}_t.
\end{equation}

The dynamics of the coalesced particle of mass $M'$ at $Z_t$ for $t \geq \tau$ are therefore
\begin{equation}
	dZ_t = -\chi m \frac{\partial V}{\partial x} \left(Z_t, x_t \right) dt + \sqrt{\frac{2\mu}{M'}} dW^{(cm)}_t,
\end{equation}
which in the presence of additional particles generalize to \eqref{eq:particle_dynamics}.
\subsection{\label{sec:subsystem_second_moment}Evolution of a subsystem's second moment}
\begin{figure}
	\captionsetup{width=\linewidth}
	\centering
	\begin{subfigure}[t]{0.49\textwidth}
		\includegraphics[width=\textwidth]{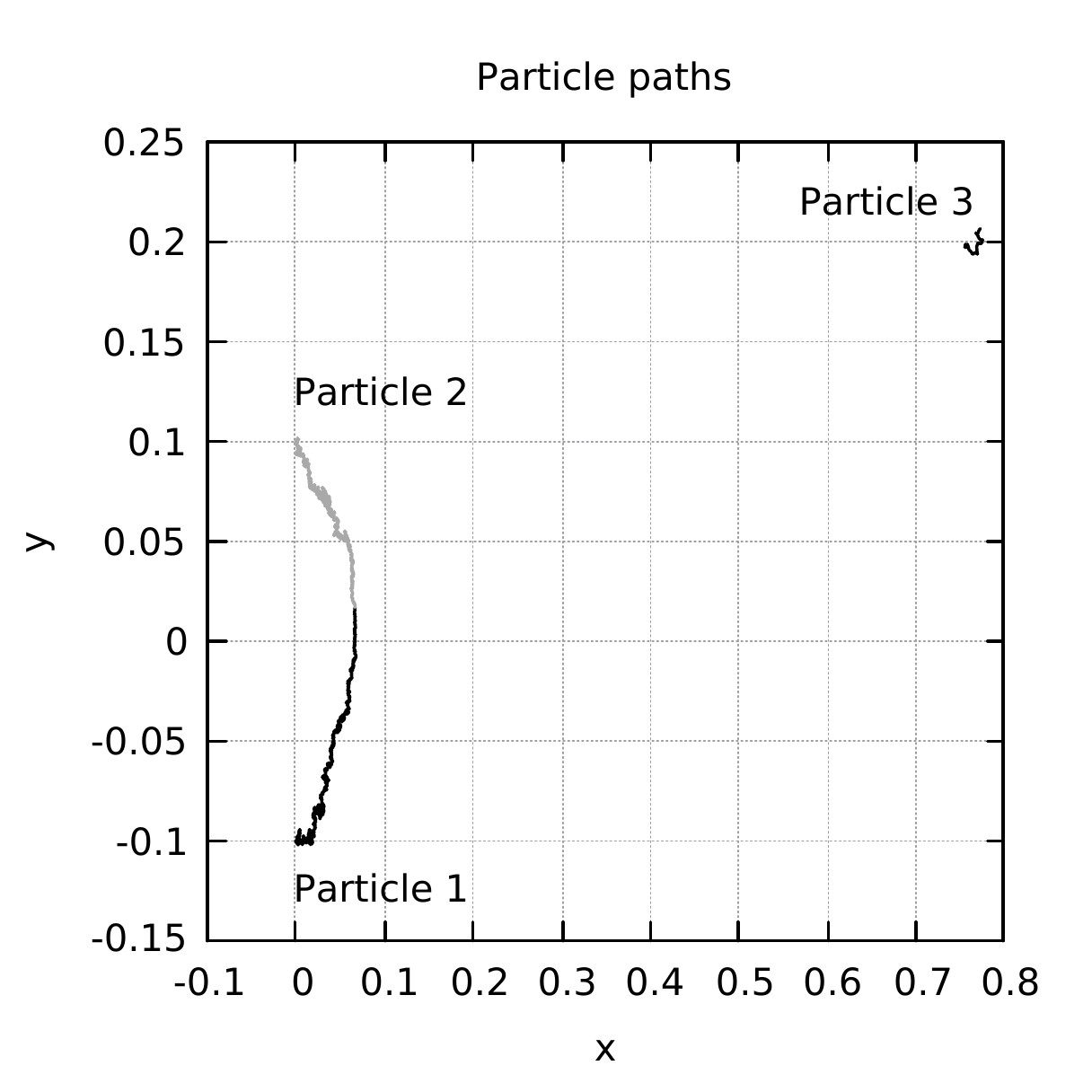}
		\caption{All three particles are drifting in the direction of the center of mass. As can be seen by the asymmetry in the paths of particles 1 and 2, the effect of the third particle on the dynamics of the first two is non-negligible.}
	\end{subfigure}
	\,
	\begin{subfigure}[t]{0.49\textwidth}
		\includegraphics[width=\textwidth]{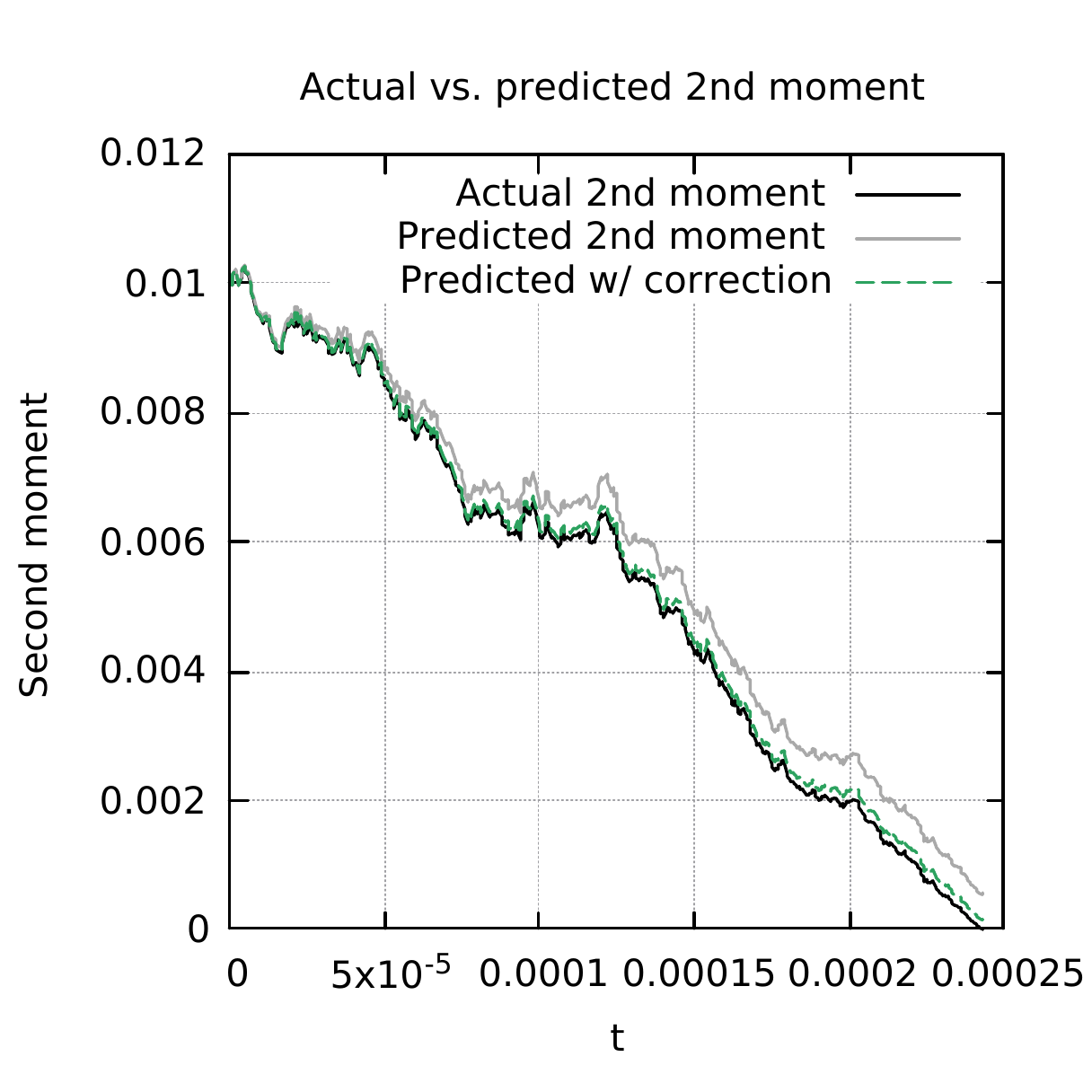}
		\caption{The second moment of the subsystem consisting of the first two particles is approximated using \eqref{eq:second_moment_zeroth} and \eqref{eq:subsystem_Y_t_monopole}, and both are compared with the real second moment. Approximation \eqref{eq:subsystem_Y_t_monopole} shows the better agreement with the actual dynamics.}
	\end{subfigure}
	\caption{An adaptive time step is used to simulate a three-particle system with $\chi~=~10$, $\tilde{\mu}~=~10$, and particle masses $m_1 = m_2 = 20$, $m_3 = 100$. The first two particles are initialized at $\left(0, \pm \frac{1}{10}\right)$, the third at $\left(\frac{4}{5} \cos \theta, \frac{4}{5} \sin \theta \right)$ with $\theta = \pi/12$.}
	\label{fig:approximated_2nd_moment}
\end{figure}
Let us compute the local second moment of the highly localized  subsystem of the first $N'$ particles in \figref{fig:formed_cluster}. First, we ignore all interactions with the outside particles not in the colliding cluster, and therefore approximate that the local second moment,
\begin{equation}
	\label{eq:second_moment_zeroth}
	\tilde{Y_t} =  \frac{1}{2(M')^2} \sum_{i,j=1}^{N'} m_i m_j \left| X^{(i)}_t - X^{(j)}_t \right|^2,
\end{equation}
evolves according to \eqref{eq:dY_t} with the summation being taken over the indices of the particle participating in the collision,
\begin{align}
	\label{eq:local_second_moment_zeroth}
	d\tilde{Y_t} \approx dQ_t = \left(\frac{4\tilde{\mu}(N'-1)}{M'} - \frac{\chi M'}{2\pi} \left(1 - \sum_{j=1}^{N'} \left( \frac{m_j}{M'} \right)^2 \right) \right) dt + 2\sqrt{\tilde{Y_t}} \sqrt{\frac{2\tilde{\mu}}{M'}} d\tilde{W}_t,
\end{align}
where
\begin{equation}
	\label{eq:weighted_subsystem_noise}
	d\tilde{W}_t = \frac{1}{(M')^{3/2} \sqrt{\tilde{Y}_t}} \sum_{i,j=1}^{N'} m_j \sqrt{m_i} \left( X^{(i)}_t - X^{(j)}_t \right) \cdot dW^{(i)}_t.
\end{equation}

As shown in \figref{fig:approximated_2nd_moment}, such an approach appears to be qualitatively correct, but introduces an error which appears to grow in time. Let us now find a higher order approximation.

As a model for the system in \figref{fig:formed_cluster}, consider a system consisting of two nearby particles of masses $m_1$ and $m_2$, and a third, distant particle of mass $m_3$, i.e. $\left| X^{(1)}_t - X^{(2)}_t \right| \ll \left| X^{(1)}_t - X^{(3)}_t \right| \approx \left| X^{(2)}_t - X^{(3)}_t \right|$. We wish to investigate how the third particle affects the second moment of the subsystem consisting of the first two particles,
\begin{equation}
	\tilde{Y_t} = \frac{m_1 m_2}{(m_1 + m_2)^2} \left| X^{(1)}_t - X^{(2)}_t \right|^2.
\end{equation}
Using \eqref{eq:particle_dynamics} and an application of Ito's lemma, we can get an exact correction to the deterministic part of the approximating process $Q_t$ given in \eqref{eq:local_second_moment_zeroth}:
\begin{align}
	\label{eq:subsystem_Y_t_exact}
	d\tilde{Y_t} = dQ_t &+ \frac{2m_1 m_2}{(m_1 + m_2)^2} \left(X^{(1)}_t - X^{(2)}_t \right) \cdot \\
	\nonumber
	& \cdot \left[-\frac{\chi m_3}{2\pi} \left(\frac{X^{(1)}_t - X^{(3)}_t}{\left| X^{(1)}_t - X^{(3)}_t \right|^2} + \frac{X^{(2)}_t - X^{(3)}_t}{\left| X^{(2)}_t - X^{(3)}_t \right|^2} \right)\right]dt.
\end{align}
We introduce the small parameter
\begin{equation}
	\epsilon_t = (X^{(2)}_t - X^{(1)}_t)/(m_1 + m_2),
\end{equation}
through which \eqref{eq:subsystem_Y_t_exact} may be approximated as
\begin{equation}
	\label{eq:approximate_second_moment}
	d\tilde{Y_t} = dQ_t - \frac{\chi m_3}{\pi} \frac{\tilde{Y}_t}{\left| X^{(cm)}_t - X^{(3)}_t \right|^2} \cos 2 \theta dt + \mathcal{O}(|\epsilon_t|^2)dt
\end{equation}
where we assume $X^{(1)}_t - X^{(3)}_t \approx X^{(2)}_t - X^{(3)}_t \approx X^{(cm)}_t - X^{(3)}_t$ and $\theta$ is the angle between $X^{(2)}_t - X^{(1)}_t$ and $X^{(cm)}_t - X^{(3)}_t$.

A similar monopole approximation may be used when there are $N-2$ particles affecting the evolution of the second moment of the first two particles. Then, 
\begin{align}
	d\tilde{Y}_t &= dQ_t - \frac{\chi \tilde{Y}_t}{\pi} \sum_{i=3}^{K+2} \frac{m_i}{\left| X^{(cm)}_t - X^{(i)}_t \right|^2} \cos 2 \theta_i + \mathcal{O}(|\epsilon_t|^2)\\
	&= dQ_t + 2\chi\tilde{Y}_t \sum_{i=3}^{K+2} m_i V''\left(\left| X^{(cm)}_t - X^{(i)}_t \right|\right) \cos 2 \theta_i + \mathcal{O}(|\epsilon_t|^2),
\end{align}
where $\theta_i$ is the angle between $X^{(2)}_t - X^{(1)}_t$ and $X^{(cm)}_t - X^{(i)}_t$, and the shorthand $V(x,y) = V(|x-y|)$ is used to simplify the expression.

By a similar argument, for an $N$ particle system with a cluster consisting of the first $N'$ particles, we have
\begin{align}
	\label{eq:subsystem_Y_t_monopole}
	d\tilde{Y_t} &\approx dQ_t + 2\chi \tilde{Y_t} \sum_{i=N'+1}^N  \sum_{j,k=1}^{N'} m_i V''\left( \left| X^{(cm)}_t - X^{(i)}_t \right| \right) \cos 2\theta_{ijk} dt, 
\end{align}
where $\theta_{ijk}$ is the angle between $X^{(j)}_t - X^{(k)}_t$ and $\tilde{X}^{(cm)}_t - X^{(i)}_t$, with
\begin{equation}
	\tilde{X}^{(cm)}_t = (m_i X^{(i)}_t + m_j X^{(j)}_t)/(m_i+m_j).
\end{equation}

Heuristically, we see that as $\tilde{Y_t} \rightarrow 0$, the corrections in \eqref{eq:subsystem_Y_t_monopole} vanish, the subsystem essentially becomes decoupled from the rest of the system, and the subsystem's second moment $\tilde{Y_t}$ becomes a squared Bessel process of negative index by \eqref{eq:negative_index}. Since the collision process (before the collision time) does not involve the creation or annihilation of particles, it appears that a highly-localized aggregate which is not decoupled from the rest of the system, but is nontheless undergoing collision, should still satisfy \eqref{eq:negative_index}, i.e.
\begin{equation}
	\label{eq:subsystem_negative_index}
	\nu(m_1, m_2, \ldots, m_{N'}) < 0,
\end{equation}
where $\nu$ is as in \eqref{eq:Bessel_index}. This informal argument suggests that for a very tight cluster, this is a sufficient condition for an aggregate to undergo collision. For a less tight cluster (even if it is separated), the contributions of the higher order corrections may prevent a collision from occurring.

\section{Simulation of particle coalescence and dynamics}
\subsection{Overview} We employ a grid-particle approach for computing interparticle interactions, which avoids pairwise computations in \eqref{eq:particle_dynamics} by introducing a continuous global potential which varies in time. We remark that similar ideas have been developed in the particle-in-cell literature (e.g. \cite{dawson1983particle}, \cite{verboncoeur2005particle}), but without coalescing stochastic particles.

We sidestep the challenge of numerically detecting singular point collisions by introducing an adaptive grid which identifies highly localized aggregates, the second moment of which is computed and simulated using the appropriate Wiener process (given by \eqref{eq:Second_moment_driving_process}) in order to identify a collision inside a timestep.

\subsection{Full numerical method}
The numerical method for the simulation of the coalescing particle system \eqref{eq:particle_dynamics} combines the upcoming sections at every timestep in the following order:
\begin{enumerate}
	\item Detect highly isolated clusters of particles with negative indices, which may collide with high probability within the upcoming time step. For each such cluster, compute the local second moment, $\tilde{Y_t}$.
	\item Simulate the particle dynamics, using adaptive timestepping when appropriate. For each particle in the above clusters, record the total increment of the driving Wiener process over the full time step.
	\item For each cluster, simulate the second moment over a time step, using \eqref{eq:Second_moment_driving_process}. If the second moment hits zero, coalesce the cluster's particles at their center of mass.
\end{enumerate}

\subsection{\label{sec:aggregate_detection}Detection of isolated aggregates}
To detect particle collisions, we first apply a density-based clustering algorithm for finding isolated particle aggregates. Such clusters are then checked for collisions, as described in \secref{sec:collision_detection}.

To find clusters, we form a coarse mesh which covers all the particles (in practice, we use a $1 \times 1$ mesh). For each cell, we compute the square of its diagonal, $s^2$, and the second moment of the particles inside the cell, $\tilde{Y}$. We call a cell ``separated'' if
\begin{equation}
	\label{eq:second_moment_separation_eta}
	\tilde{Y}/s^2 < \eta \ll 1,
\end{equation} 
where $\eta$ is some fixed constant (in practice, the authors use $\eta = 0.1$). If a cell is not separated, and has more than two particles, then we refine the cell into four equally-sized cells, and repeat this procedure with each subcell.

A separated cell is kept if it is ``collidable,'' otherwise it is refined as well. A cell is collidable if its index $\nu$ is negative, and the second moment satisfies
\begin{equation}
	\tilde{Y} + \alpha \Delta t + 2\beta\sqrt{\tilde{Y}}\Phi^{-1}(p)\sqrt{\Delta t} < 0,
\end{equation}
where $\alpha$ and $\beta$ are given as in \eqref{eq:alpha_beta}, $\Phi$ is the normal distribution function, and $0 < p \ll 1$ is some small probability. The interpretation of this inequality is that it excludes cells which may collide within the time step with very low probabilities.




\subsection{\label{sec:simulating_dynamics}Particle dynamics}
\begin{figure}[t]
	\centering
	\includegraphics[scale=0.5]{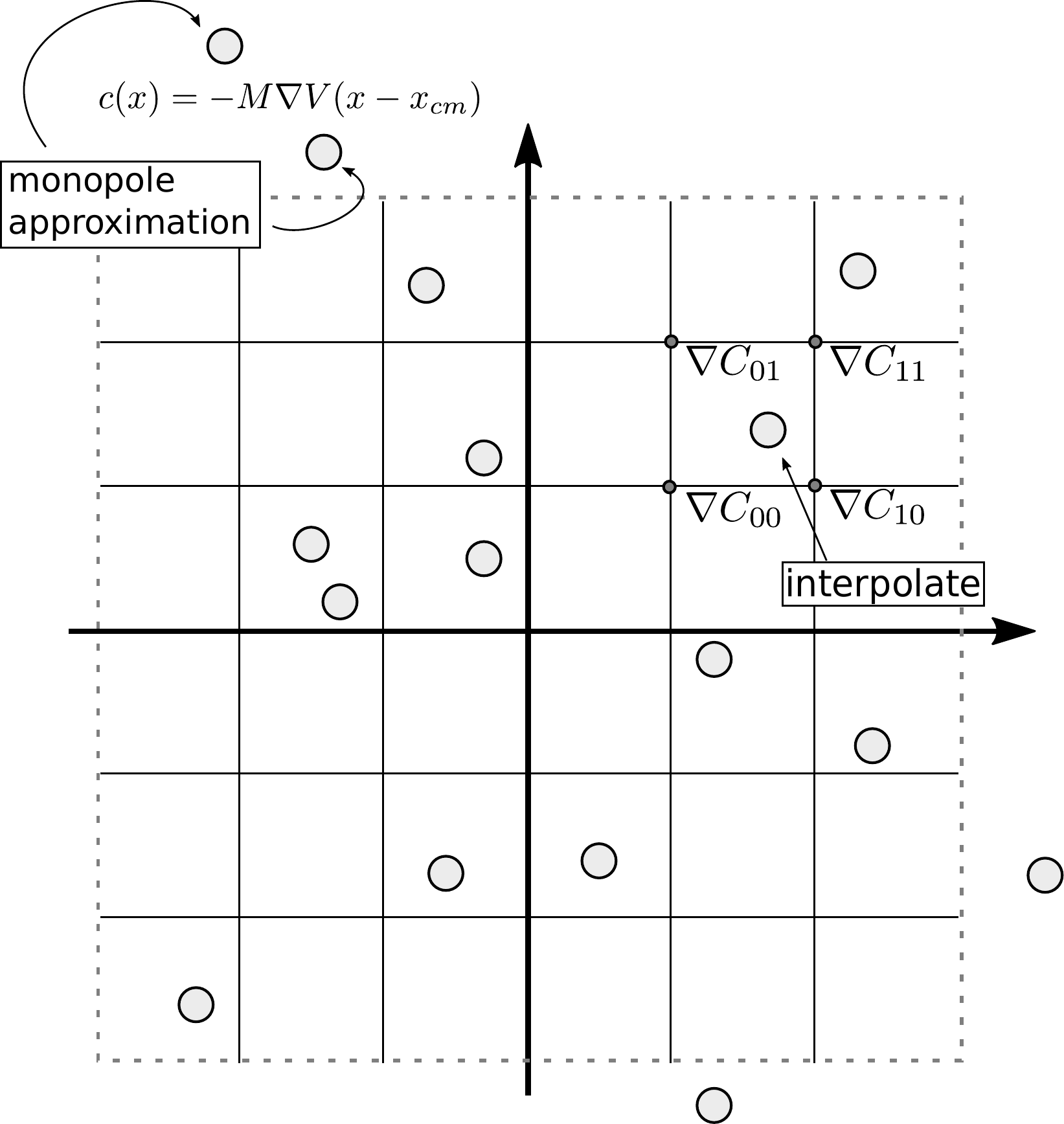}
	\caption{Inside the computational domain, which we denote by the dashed box, $\nabla C_{ij}$ is computed numerically, and then bilinearly interpolated at the point inside the cell. Outside the computational domain, we approximate $\nabla c$ via a monopole approximation ($x_{cm}$ denotes the center of mass).}
	\label{fig:c_definition}
\end{figure}
Since $V(x,y) = \frac{1}{2\pi}\ln|x-y|$ is the fundamental solution to the Laplace operator, we can get a global potential for the the dynamics given in \eqref{eq:particle_dynamics}:
\begin{align}
	\begin{cases}
	dX^{(n)}_t &= \chi \nabla c\left(X^{(n)}_t, t\right) dt + \sqrt{\frac{2\tilde{\mu}}{m_n}} dW^{(n)}_t, \\
	\Delta c &= -P, \\
	P(x,t) &= \sum_{i=1}^N m_i \delta \left( x - X^{(i)}_t \right),
	\end{cases}
\end{align}
where $c$ is the interaction  field and $P$ is the empirical mass density. The case of a different $V$ can be treated similarly.

For the simulation of these particle dynamics, we discretize a computational domain as in figure  \ref{fig:c_definition}, and use the particles to interpolate a mass density field $P_{ij}$ onto the field. We then numerically solve for the mean field $C_{ij}$. To advance by $\Delta t$ forward in time, we introduce adaptive time steps $\Delta \tau_1, \ldots, \Delta \tau_{K(n,t)}$ (this is needed for stability reasons---see below), and use a forward Euler-Maruyama scheme to simulate the dynamics of each particle:
\begin{equation}
	 X^{(n)}(t + \Delta \tau_i) =  X^{(n)}(t) + \chi \nabla c\left(X^{(n)}(t), t\right) \Delta \tau_i + \sqrt{\frac{2\tilde{\mu}}{m_n}}N_i^{(n)}(0,1) \sqrt{\Delta \tau_i},
\end{equation}
where $N_i^{(n)}(0,1)$ is a normal Gaussian random variable. This bookkeeping of the noise is helpful for numerically detecting collisions, where we need the quantity
\begin{equation}
	\label{eq:nth_increment}
	\Delta W^{(n)}(t) = \sum_{i=1}^{K(n,t)} \sqrt{\Delta \tau_i} N^{(n)}_i(0,1),
\end{equation}
i.e. the increment of the $n$th Wiener process $W^{(n)}_t$ between $t$ and $t+\Delta t$ (see \secref{sec:collision_detection}).

We approximate $\nabla c(x)$ in two steps. First, we construct the gradient field $\nabla C_{ij}~=~(CX_{ij}, CY_{ij})$ using the second order approximation
\begin{align}
	CX_{ij} &= \frac{C_{i+1,j} - C_{i-1,j}}{2\Delta x}, \\
	CY_{ij} &= \frac{C_{i,j+1} - C_{i,j-1}}{2\Delta x}.
\end{align}
Then we approximate $\nabla c(x,t)$ using a bilinear interpolation of the values of $\nabla C$ at the four nearest grid points. In the case that $x$ is not inside the computational domain, we use a monopole approximation:
\begin{equation}
	\nabla c(x) = -M \nabla V(x-x_{cm}),
\end{equation}
where $x_{cm}$ is the center of mass of the system. Since the primary novelty of this numerical method is in its applicability to colliding systems, an appropriately-chosen computational domain (i.e. one which overlaps with most of the mass of the system) will make use of the monopole approximation rarely. Nonuniform meshes may be used as well, but have not been observed to make a significant improvement in systems with most of the mass sufficiently away from the boundaries .

As described in \cite{fatkullin2013study}, an adaptive time step is dynamically chosen such that the expected length of a particle's jump does not exceed the mesh size. This is necessary to prevent spurious mass oscillations around singularities in $c(x,t)$. Since $V$ has logarithmic singularities, we expect that time steps can get as small as $\Delta \tau \sim \Delta x^2$.

The mass density field $P_{ij}$ is computed by bilinearly interpolating the mass of each particle onto the four nearest grid points. The result is divided by $\Delta x^2$, to get a mass density. This first-moment-preserving approach prevents particles from ``self-interacting,'' a phenomenon which creates an artificial flux towards grid points, as described in \cite{fatkullin2013study}.

The mean field $C_{ij}$ is solved for on the computational domain using a standard finite-differences scheme:
\begin{equation}
	\frac{1}{\Delta x^2} \left( C_{i+1,j} + C_{i-1,j} + C_{i,j-1} + C_{i,j+1} - 4C_{ij} \right) = -P_{ij}.
\end{equation}
The monopole approximation is used for the boundary conditions:
\begin{equation}
	C_{ij} = -MV(X_{ij} - x_{cm})
\end{equation}
for $X_{ij}$ on the boundary on the computational domain, and $x_{cm}$ the center of mass of the particle system.

\subsection{\label{sec:collision_detection}Detection of collisions in isolated aggregates} After all the particles are propagated over one time step, we consider the terminal cells returned by the algorithm given in \secref{sec:aggregate_detection}. We approximate the evolution of the second moment inside each cell which is both separated and collidable. To do this, for each cell, we compute the quantity
\begin{equation}
	\Delta \tilde{Y} = \alpha \Delta t + 2 \beta \sqrt{\tilde{Y}} \Delta \tilde{W}_t,
\end{equation}
and coalesce all the particles at their new center of mass if $\Delta \tilde{Y} \leq 0$.

The increment $\Delta \tilde{W}_t$ is given by \eqref{eq:weighted_subsystem_noise}, i.e.
\begin{equation}
	\label{eq:numerical_driving_process_Y_t}
	\Delta \tilde{W}_t = \frac{1}{(M')^{3/2} \sqrt{\tilde{Y}_t}} \sum_{i,j=1}^{N'} m_j \sqrt{m_i} \left( X^{(i)}_t - X^{(j)}_t \right) \cdot \Delta W^{(i)}(t).
\end{equation}
The cost of computing the above sum can be significantly reduced using the following identity:
\begin{align}
	dW_t &= \frac{1}{\sqrt{MY_t}} \sum_{i=1}^{N} \sqrt{m_i} \left( X^{(i)}_t - X^{(cm)}_t \right) \cdot dW^{(i)}_t,
\end{align}
from which
\begin{equation}
	\label{eq:simplied_weighted_noise}
	\Delta \tilde{W}_t = \frac{1}{\sqrt{M \tilde{Y}_t}} \sum_{i=1}^{N'} \sqrt{m_i} \left( X^{(i)}_t - X^{(cm)}_t \right) \cdot N^{(i)}(t) \sqrt{\Delta t},
\end{equation}
where $X^{(cm)}_t$ is the center of mass of the cell.

We note the dynamics of the second moment may be approximated more accurately by taking advantage of the first order correction presented in \secref{sec:subsystem_second_moment}, but the necessity of such corrections may be avoided by simply choosing a very small localization parameter $\eta$, as in \eqref{eq:second_moment_separation_eta}.

\section{Finite-time blow-up in hydrodynamic limits}
\subsection{Overview}
We first show how the PKS particle system described in the introduction fits in the context of the present work. We then formally derive the hydrodynamic limit of a particle system with masses approaching zero nonuniformly, which we call the multispecies Patlak-Keller-Segel system (MPKS), and derive a finite-time blow-up condition. Finally, we show how the hydrodynamic limit of the system may be taken in such a way that the limit is a regularized MPKS system after the time of blow-up.

\subsection{\label{sec:PKS}The Patak-Keller-Segel particle system}
As already described in \secref{sec:PKS_intro}, the PKS is given by the following system of PDEs:
\begin{align}
\label{eq:PKS}
\begin{cases}
	\partial_t \rho &= \nabla \cdot ( \mu \nabla \rho - \chi \rho \nabla c ), \\
	\Delta c &= -\rho,
\end{cases}
\end{align}
where the boundary condition $\rho(x,t) \rightarrow 0$ as $|x| \rightarrow \infty$ is enforced, and mass is conserved: $\int \rho(x,t)dx = M$.

The PKS may be rewritten more compactly as an integrodifferential equation:
\begin{equation}
	\label{eq:PKS_integrodifferential}
	\partial_t \rho = \nabla \cdot (\mu \nabla \rho + \chi \rho \nabla (V \ast \rho)),
\end{equation}
where $V(x) = \frac{1}{2\pi}\ln|x|$, as before, is the fundamental solution of the Laplace operator. Observe that if $c$ is predetermined, then the first equation in \eqref{eq:PKS} is the Fokker-Planck equation for the process
\begin{equation}
	\label{eq:PKS_SDE}
	dX_t = \chi \nabla c(X_t, t) dt + \sqrt{2\mu} dW_t.
\end{equation}
It follows that for an $N$-particle system with positions $X^{(n)}_t$, the empirical mass density
\begin{equation}
	P_N (x,t) = \frac{M}{N} \sum_{n=1}^N \delta \left(x - X^{(n)}_t \right)
\end{equation}
approximates the solution to the PKS $\rho$. 

Since $\nabla c$ is unknown, we approximate it by the mean field created by the particles themselves: this is readily done making the substitution $c \rightarrow -V \ast P_N$, as suggested by \eqref{eq:PKS_integrodifferential}. We arrive at
\begin{equation}
	\label{eq:PKS_particle}
	dX^{(n)}_t = -\frac{\chi M}{N} \frac{\partial}{\partial X^{(n)}_t} \sum_{i \neq n} V(X^{(n)}_t, X^{(i)}_t) dt + \sqrt{2\mu} dW^{(n)}_t.
\end{equation}
This is simply the particle system described in the bulk of this work, with $m_n = M/N$ and the diffusion coefficient 
\begin{equation}
	\label{eq:PKS_diffusion_coefficient}
\tilde{\mu} = \mu m_n = \frac{\mu M}{N}.
\end{equation}
Thus, the PKS with total mass $M$ and diffusion coefficient $\mu$ can be viewed as the hydrodynamic limit of the particle system with the above parameters. 

The particle system described in this work collides only when the index of the system \eqref{eq:Bessel_index} is negative. Similarly, the PKS forms singularities when the total mass is above the critical mass $M_c = 8\pi\mu/\chi$ \cite{blanchet2006two}. Let us show that these two criteria coincide in the hydrodynamic limit.

Substituting the necessary diffusion coefficient \eqref{eq:PKS_diffusion_coefficient} into the definition of the Bessel index \eqref{eq:Bessel_index}, we get the PKS index:
\begin{align}
	\nu_{PKS} &= (N-2) - \frac{\chi M N}{8\pi\mu} \left(1 - \sum_k \left( \frac{m_k}{M} \right)^2 \right) \\
	\label{eq:MPKS_index}
	&= N \left[ \left(1 - \frac{2}{N} \right) - \frac{\chi M}{8\pi\mu} \left(1 - \sum_k \left( \frac{m_k}{M} \right)^2 \right) \right] \\
	\label{eq:PKS_index}	
	&= N \left[ \left(1 - \frac{2}{N} \right) - \frac{\chi M}{8\pi\mu} \left(1 - \frac{1}{N}  \right) \right] \\
	\label{eq:PKS_index_convenient}	
	&= (N-1)\left(1 - \frac{\chi M}{8\pi\mu}\right) - 1.
\end{align}
As per the classification of the origin for the second moment, listed in \secref{sec:collision_criterion}, we have that a finite-time collision will occur when $\nu \leq -1$. This criterion applied to \eqref{eq:PKS_index_convenient} reduces exactly to $M > 8\pi\mu/\chi$---the necessary and sufficient condition for finite-time blow-up in the PKS.
 
\subsection{Post-blow-up PKS and particle coalescence} The PKS has been regularized and investigated post-blow-up in several works, including \cite{velazquez2004pointsI, velazquez2004pointsII} and \cite{dolbeault2009two}. Although the post-blow-up dynamics are slightly different in the two works, they share the common feature that the density becomes a measure, and splits into a regular, and an atomic component consisting of $K_t$ point masses:
\begin{equation}
	\rho(x,t) = \rho_{reg}(x,t) + \sum_{n=1}^{K_t} M_n(t) \delta \left( x-x^{(n)}_t \right),
\end{equation}
where the $n$th atomic component has a smoothly-evolving mass $M_n(t) \geq 8\pi\mu/\chi$, supported on a point moving along a smooth path. The point masses may emerge or collide, and thus their number $K_t$ varies in time. Mass is locally transferred from the regular component to each atomic component as
\begin{equation}
	\label{eq:mass_transfer}
	\frac{dM_n}{dt} = \rho_{reg}\left(x^{(n)}_t,t\right) M_n.
\end{equation}
With these dynamics, it can be shown \cite{dolbeault2009two} that the second moment of this system evolves as
\begin{equation}
	\label{eq:regularized_second_moment}
	\frac{d}{dt} \left( \frac{1}{M} \int |x|^2 \rho(x,t) dx \right) = 4\mu\frac{\bar{M}}{M} - \frac{\chi M}{2\pi} \left( 1 - \sum_{i=1}^{K_t} \left( \frac{M_i(t)}{M} \right)^2 \right),
\end{equation}
where $\bar{M} = M - \sum_{i=1}^{K_t} M_i(t)$ is the mass of the regular component (we note the quantity of interest in the PKS literature is typically the unnormalized second moment, which we choose to normalize, due to its geometric interpretation).

In the context of the PKS particle system, we expect light, uncoalesced particles to correspond to the regular component of the solution to the PKS, and each massive, coalesced particles to correspond to point mass in the atomic component of the solution to the PKS. By the previous section, such particles should only have mass above $8\pi\mu/\chi$, as in the PKS. Let us recover equation \eqref{eq:regularized_second_moment} using the particle system, assuming that this correspondence is true.

Consider a PKS system with smooth initial conditions, which blows up in finite time, and has an atomic component of mass $M_1$, consisting of one point mass, at time $t=T$. Now consider a PKS particle system, initialized with $N_0$ particles distributed according to the initial conditions given to the PKS PDE. The second moment $Y_t$ evolves according to
\begin{equation}
	dY_t = \alpha dt + 2 \beta \sqrt{Y_t} dW_t,
\end{equation}
where $\alpha$ and $\beta$ are given in \eqref{eq:alpha_beta}, with $\tilde{\mu} = \mu M / N_0$, and $m_n = M/N_0$ initially. Near $t=T$, there should be one massive particle, consisting of $k$ coalesced light particles. Plugging this into \eqref{eq:alpha_beta}, we get:
\begin{align}
	\alpha &= \frac{4\mu}{M} \frac{M}{N_0}\left(N_0 - k + 1 - 1\right) - \frac{\chi M}{2\pi} \left( 1 - \frac{N_0 - k}{N_0^2} - \left(\frac{M_1}{M}\right)^2\right) \\
	&= \frac{4\mu}{M} \left( \frac{N_0 - k}{N_0} M \right) - \frac{\chi M}{2\pi} \left( 1 - \frac{N_0 - k}{N_0^2} - \left(\frac{M_1}{M}\right)^2\right).
\end{align}
As $N_0 \rightarrow \infty$, we get $\beta \rightarrow 0$, and $Y_t$ becomes deterministic:
\begin{equation}
	dY_t \rightarrow 4\mu \frac{\bar{M}}{M} - \frac{\chi M}{2\pi} \left(1 - \left( \frac{M_1}{M} \right)^2 \right) dt,
\end{equation}
consistent with \eqref{eq:regularized_second_moment} for a single point mass. A similar argument can be used to derive \eqref{eq:regularized_second_moment} fully.

\subsection{\label{sec:MPKS_limit}Hydrodynamic limit to the multispecies PKS model}
We remark that the sign of the PKS particle system's index \eqref{eq:PKS_index}  becomes independent of $N$ as $N \rightarrow \infty$. This convenient property occurs only because $\tilde{\mu} \sim 1/N$, and is actually independent of the the particle masses, as long as the total sum of the particle masses is fixed and the mass of each individual particle approaches zero.  Thus the question of the limiting system when individual particles approach $0$ nonuniformly arises naturally.

As a first basic example, let us consider the system
\begin{equation}
	dX^{(n)}_t = - \chi \frac{\partial}{\partial X^{(n)}_t} \sum_{i \neq n} m_i V(X^{(n)}_t, X^{(i)}_t) dt + \sqrt{\frac{2\mu M}{N m_n}} dW^{(n)}_t,
\end{equation}
where $N=2N'$, $M=M_1+M_2$, $m_i~=~M_1/N'$ for $i \leq N'$ and $m_i~=~M_2/N'$ for $i > N'$. That is, we break up the system into two families, the first family containing $N'$ particles of uniform mass $m_a = M_1/N'$, and the second family containing $N'$ particles of uniform mass $m_b = M_2/N'$. The particle dynamics are then given by
\begin{align}
	\begin{cases}
	dX^{(n)}_t &= \chi \nabla c(X^{(n)}_t, t) dt + \sqrt{\mu\left(1+ \frac{M_1}{M_2} \right)}dW^{(n)}_t, \,\,\,\, n \leq N' \\
	dX^{(n)}_t &= \chi \nabla c(X^{(n)}_t, t) dt + \sqrt{\mu\left(1+\frac{M_2}{M_1} \right)}dW^{(n)}_t, \,\,\,\, n > N' \\
	\Delta c &= -P_1(x) - P_2(x),
	\end{cases}
\end{align}
where $P_1$ and $P_2$ are the empirical mass densities of the particles of the first and second mass:
\begin{align}
	P_1(x) &= \sum_{i=1}^{N'}m_a \delta \left(x-X^{(i)}_t \right), \\
	P_2(x) &= \sum_{i=N'+1}^{2N'}m_b \delta \left(x-X^{(i)}_t \right).
\end{align}
Appealing once more to the formal derivation of the hydrodynamic limit described earlier, we expect that $P_i$ approximates $\rho_i$ in the limit $N \rightarrow \infty$, where
\begin{equation}
	\begin{cases}
		\partial_t \rho_1 &= \nabla \cdot \left( \frac{\mu}{2} \left(1+\frac{M_2}{M_1} \right) \nabla \rho_1 - \chi \rho_1 \nabla c \right), \\
		\partial_t \rho_2 &= \nabla \cdot \left( \frac{\mu}{2} \left(1+\frac{M_1}{M_2} \right) \nabla \rho_2 - \chi \rho_2 \nabla c \right), \\
		\Delta c &= -(\rho_1 + \rho_2).
	\end{cases}
\end{equation}
The above system can be seen as a ``two species'' PKS model, in which two species attract each other through the same mechanism, but have different average diffusion rates.

Similarly, we may break the system up into $K$ families, each family of total mass $M_i$ and containing $N_i$ particles of uniform mass $M_i/N_i$. We take the hydrodynamic limit by fixing $\eta_i > 0$ for $1 \leq i \leq K$ such that
\begin{equation}
	\label{eq:hydrodynamic_mass_ratios1}
	\eta_1 + \cdots + \eta_K = 1,
\end{equation}
and letting $N \rightarrow \infty$ in such a way that
\begin{equation}
	\label{eq:hydrodynamic_particle_numbers}
	N_i = \eta_i N.
\end{equation}
Then $P_i \rightarrow \rho_i$, where
\begin{align}
	\label{eq:nonuniform_masses_hydrodynamic}
	\begin{cases}
		\partial_t \rho_1 &= \nabla \cdot \left( \mu_1 \nabla \rho_1 - \chi \rho_1 \nabla c \right), \\
		&\vdots \\
		\partial_t \rho_K &= \nabla \cdot \left( \mu_K \nabla \rho_K - \chi \rho_K \nabla c \right), \\
		\Delta c  &= -(\rho_1 + \cdots + \rho_K),
	\end{cases}
\end{align}
with $\int \rho_i = M_i$ and
\begin{equation}
	\label{eq:nonuniform_masses_diffusion_coefficients}
	\mu_i = \frac{M}{M_i} \eta_i \mu = \lim_{N\rightarrow \infty} \frac{M/N}{M_i/N_i} \mu,
\end{equation}
which can be interpreted as $\mu$ scaled by the ratio of the overall system's average particle mass, to the $i$th family's particle mass. We will refer to \eqref{eq:nonuniform_masses_hydrodynamic} as the ``multispecies Patlak-Keller-Segel system'' (MPKS).

The excluded case $\eta_i = 0$ corresponds to a mass of $M_i$ being supported entirely on a singular component of the solution post-blow-up.
\subsection{\label{sec:singularity_formation}Formation of singularities in the MPKS} 
As can be seen from \eqref{eq:MPKS_index}, the sign of the index of a particle system that's taken to its hydrodynamic limit becomes independent of the number of particles, and can therefore fully collide in finite time, if a specific mass condition is satisfied. In the PDE, this corresponds to a finite-time blow up. Let us verify that this is indeed the case.

Assume arbitrary diffusion coefficients $\mu_i$. Let $P(x,t) = \sum_{i=1}^K \rho_i(x,t)$ be the total mass density of an MPKS system. Then $\int_{\mathbb{R}^2} P(x,t) dx = M_1 + \cdots + M_K = M$, and
\begin{equation}
	c(x,t) = -\frac{1}{2\pi}\int_{\mathbb{R}^2}\ln|x-y|P(y)dy.
\end{equation}
To show the existence of finite-time blow-up, define the second moment of the system,
\begin{equation}
	F(t) = \int_{\mathbb{R}^2} P(x,t) |x|^2 dx,
\end{equation}
and compute its derivative:
\begin{align}
\label{eq:MPKS_second_moment_derivate}
	F'(t) = \sum_{i=1}^K \left(4\mu_i - \frac{\chi M}{2\pi}\right)M_i,
\end{align}
where the detailed computation is given in \appref{app:MPKS_second_moment_derivation}. Thus for constants satisfying
\begin{equation}
	\label{eq:MPKS_blow_up_condition}
	\sum_{i=1}^K \left(4\mu_i - \frac{\chi M}{2\pi}\right)M_i < 0,
\end{equation}
the second moment vanishes in finite time, but the total mass is conserved--thus implying the formation of a singularity.

As an aside, we remark that the the formula given by \eqref{eq:MPKS_second_moment_derivate} remains valid when each component has a different chemosensitivity $\chi_i$. Furthermore, we note that the blow-up condition \eqref{eq:MPKS_blow_up_condition} is satifised when $M > \max (8\pi\mu_i / \chi)$, i.e. the MPKS forms a singularity when its total system mass is greater than the classic PKS critical mass for each separate components. Recalling the special structure of the diffusion coefficients in the hydrodynamic limit of the particle system \eqref{eq:nonuniform_masses_diffusion_coefficients}, we see that the blow-up condition \eqref{eq:MPKS_blow_up_condition} coincides with the full particle system collision condition $\nu_{PKS} < 0$, where $\nu_{PKS}$ is as in \eqref{eq:MPKS_index}. 

For two species, the system was investigated in \cite{conca2011remarks}, where initial data were classified in terms of having solutions which either blow up in finite time, or are global in time. Interestingly, that work showed that there exist initial data corresponding to finite time blow-up, for which the second moment is increasing, i.e. $F'(t)>0$---in analogy with \eqref{eq:index_difference}. An optimal classification was obtained for a disc domain in \cite{espejo2013sharp}, though questions, such as if blow up occurs simultaneously in all components, remain (this question was affirmatively answered for the radial case in \cite{espejo2009simultaneous}). In \secref{sec:MPKS_experiment}, we investigate how the second moments of components of the two species MPKS evolve in the regime that a singularity forms in finite time with $F'(t)>0$.

We expect that the MPKS can be regularized past blow-up times using a singular perturbation limit, as was done in \cite{velazquez2004pointsI, velazquez2004pointsII} for the PKS, and proposed in \cite{kurganov2014numerical} for the MPKS. In this case, the presented method is well-suited for the investigation of this regularization.

\begin{figure}
	\includegraphics[width=\linewidth]{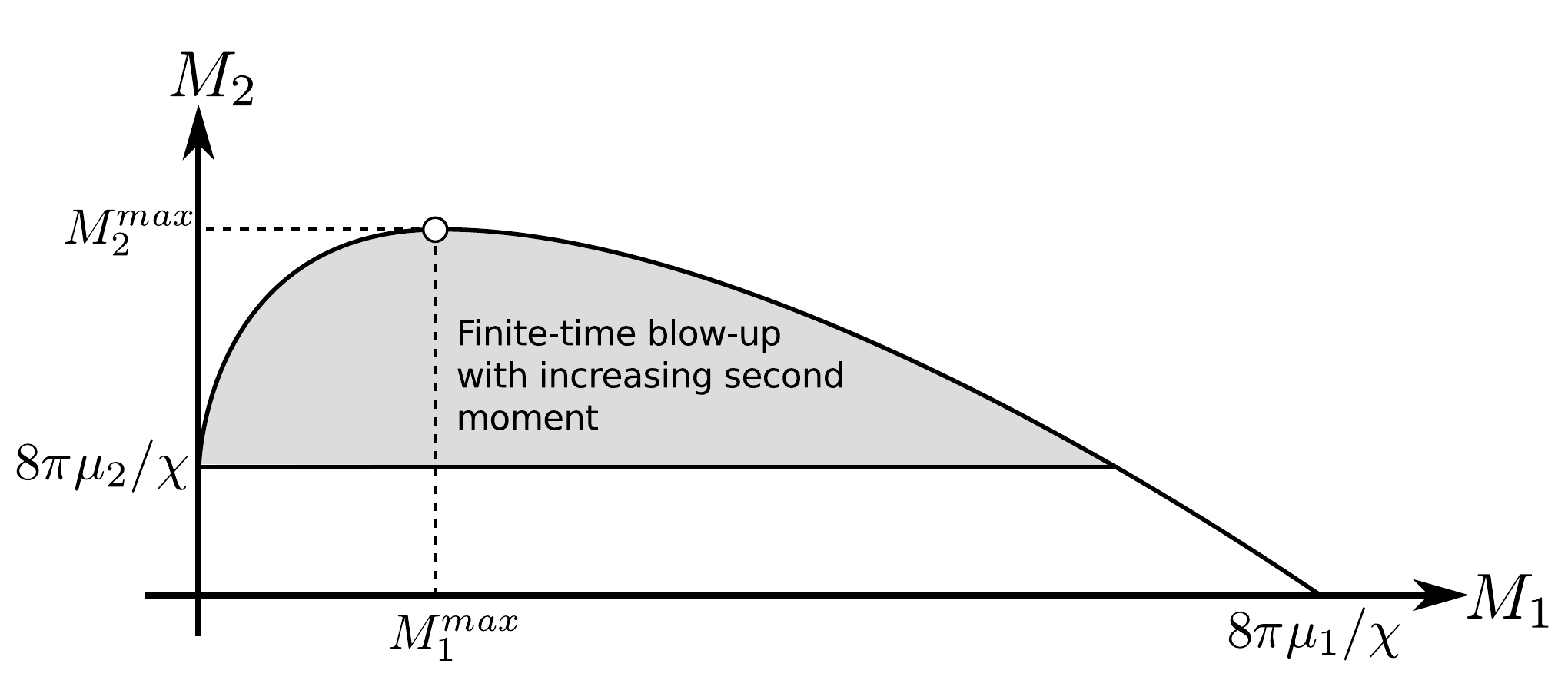}

	\caption{For the two species MPKS system, the second moment increases when the point $(M_1, M_2)$ lies below the curve obtained by setting the right hand side of \eqref{eq:MPKS_second_moment_derivate} to zero. However, it was shown in \cite{conca2011remarks} that finite-time blow-up will occur for radially-symmetric initial data when $M_2 > 8\pi\mu_2/\chi$; thus, unlike in the PKS, it is possible (when $\mu_1 > 2\mu_2$) for a system to both spread across the plane, and form a singularity in finite time. The values of $M_1^{max}$ and $M_2^{max}$ are given in \eqref{eq:MPKS_M_max}. A typical region in which this atypical behavior occurs is shaded above, with parameters $\chi = 100, \mu_1 = 10, \mu_2 = 1$. In the aforementioned work, it was hypothesized that the second moment of one component increases, while the second moment of the other component decreases. We investigate this possibility in \secref{sec:MPKS_experiment}.}
	\label{fig:MPKS_blow_up_mass}
\end{figure}

\subsection{More general $V$}
As the particle system dynamics are equally valid for choices of $V$ which are not scaled logarithms, we left some formulas somewhat general, simply in terms of the derivatives of $V$. Particle \emph{coalescence}, however, strongly depends on there being a logarithmic singularity in $V$. This is necessary to connect collisions to the Bessel process.

We note that, in the plane, the fundamental solution to a radially-symmetric, elliptic operator $\mathcal{L}$ with sufficiently regular coefficients (as in \eqref{eq:elliptic_operator}) has logarithmic singularities. It therefore follows that the discussion above applies in the case when $V$ is such a fundamental solution. That is, suppose $V(x,y) \sim \gamma \ln|x-y|$ as $|x-y|\rightarrow 0$. Then the index formulas used in the previous sections should be replaced by the following index:
\begin{equation}
	\nu^\mathcal{L}(m_1, m_2, \cdots, m_N) = N \left(1 - \frac{2}{N} \right) - \frac{\gamma \chi M^2}{4 \tilde{\mu}} \left( 1 - \sum_{j} \left(\frac{m_j}{M}\right)^2 \right).
\end{equation}

Applying the same procedure as in \secref{sec:MPKS_limit} will result in a hydrodynamic limit which solves
\begin{align}
	\label{eq:nonuniform_masses_hydrodynamic_general}
	\begin{cases}
		\partial_t \rho_1 &= \nabla \cdot \left( \mu_1 \nabla \rho_1 - \chi \rho_1 \nabla c \right), \\
		&\vdots \\
		\partial_t \rho_K &= \nabla \cdot \left( \mu_K \nabla \rho_K - \chi \rho_K \nabla c \right), \\
		\mathcal{L} c  &= -(\rho_1 + \cdots + \rho_K),
	\end{cases}
\end{align}
with post-blow-up dynamics similar to the ones given for the PKS in  \cite{velazquez2004pointsI, velazquez2004pointsII} and \cite{dolbeault2009two}.

\section{Numerical simulations\label{sec:numerical_simulations}}
\subsection{Overview}
One application of this work is in developing a numerical method for the PKS and PKS-like systems, which is able to handle the formation of singularities, as well as post-blow-up dynamics. Let us consider two example applications, for which we explicitly know the expected behavior: the evolution of the second moment for the PKS, pre- and post-blow-up, as given in \eqref{eq:regularized_second_moment}, and blow-up with an increasing second moment in the two species MPKS, as described in \figref{fig:MPKS_blow_up_mass}. In the first, we will show that the second moment of our particle approximation evolves as predicted by \cite{dolbeault2009two} both before and after blow-up, confirming that our numerical method correctly transitions from approximating smooth solutions to the PKS, to approximating measure-valued solutions. In the second, we will see how the second moment of the components of a two species MPKS system with masses inside the shaded region in \figref{fig:MPKS_blow_up_mass} evolve, thus giving numerical evidence to the idea that blow-up in this regime occurs via one contracting, and one expanding component.

We remark the presented numerical method is parallelizable, and scales approximately linearly with the number of particles. It can therefore be used to simulate a large number (on the order of millions) of particles very quickly. Averaging over such large ensembles reduces observed stochastic fluctuations to a minimum, as may be noted from the examples in this section.

\subsection{Regularized PKS\label{sec:regularized_PKS_experiment}} For the first example, we reproduce the equation \eqref{eq:regularized_second_moment} for the PKS second moment:
\begin{equation}
	\frac{d}{dt} \left( \frac{1}{M} \int |x|^2 \rho(x,t) dx \right) = 4\mu\frac{\bar{M}}{M} - \frac{\chi M}{2\pi} \left( 1 - \sum_{i=1}^{K_t} \left( \frac{M_i(t)}{M} \right)^2 \right).
\end{equation}
Thus, the graph of the second moment of a critical PKS system will initially appear linear, then decelerate, and then--depending on the mass distribution--will either become linear again (with a different slope), or continually change its slope due to nonstop mass transfer to the atomic component. Using the numerical method developed in this work, this second moment evolution can be observed. For a PKS system with mobility $\mu$ and chemosensitivity $\chi$, we associate an $N_0$-particle coalescing particle system with $\tilde{\mu} = \mu M / N_0$ and $m_n = M/N_0$, and approximate $\rho$ by the empirical mass density. As this particle approximation has been shown to be effective in approximating the PKS pre-blow-up \cite{haskovec2009stochastic, fatkullin2013study}, we specifically concentrate on the formation and detection of singularities.

\subsubsection{Mass transfer to singularity \label{sec:mass_transfer_to_singularity}} 
\begin{figure}
\begin{center}
	\includegraphics[width=0.85\linewidth]{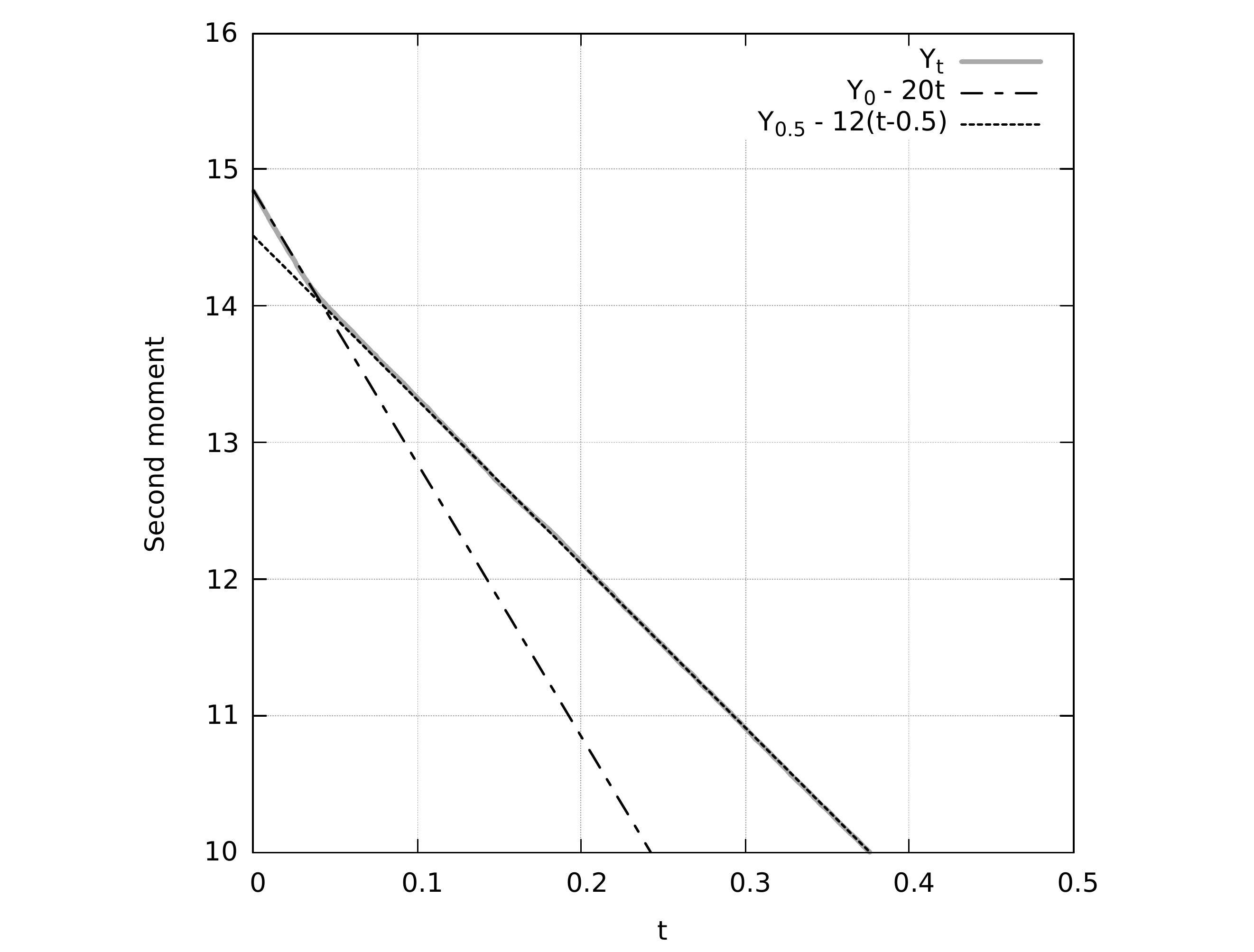}
\end{center}
	\caption{We simulate $40 \times 10^3$ particles to approximate the system described in \secref{sec:regularized_PKS_experiment}. Initially, the second moment decreases at a rate of $-20$, as predicted by the classic PKS formula for blow-up. Near $t=0.05$, a singularity is formed, and the slope of the graph of second moment suddenly changes to $-12$. Each dashed line is fitted to one only  point---i.e. the particle approximation of the PDE is effective post-blow-up.}
	\label{fig:PKS_simulation_zoom}
\end{figure}

\begin{figure}
	\begin{center}
		\includegraphics[width=\linewidth]{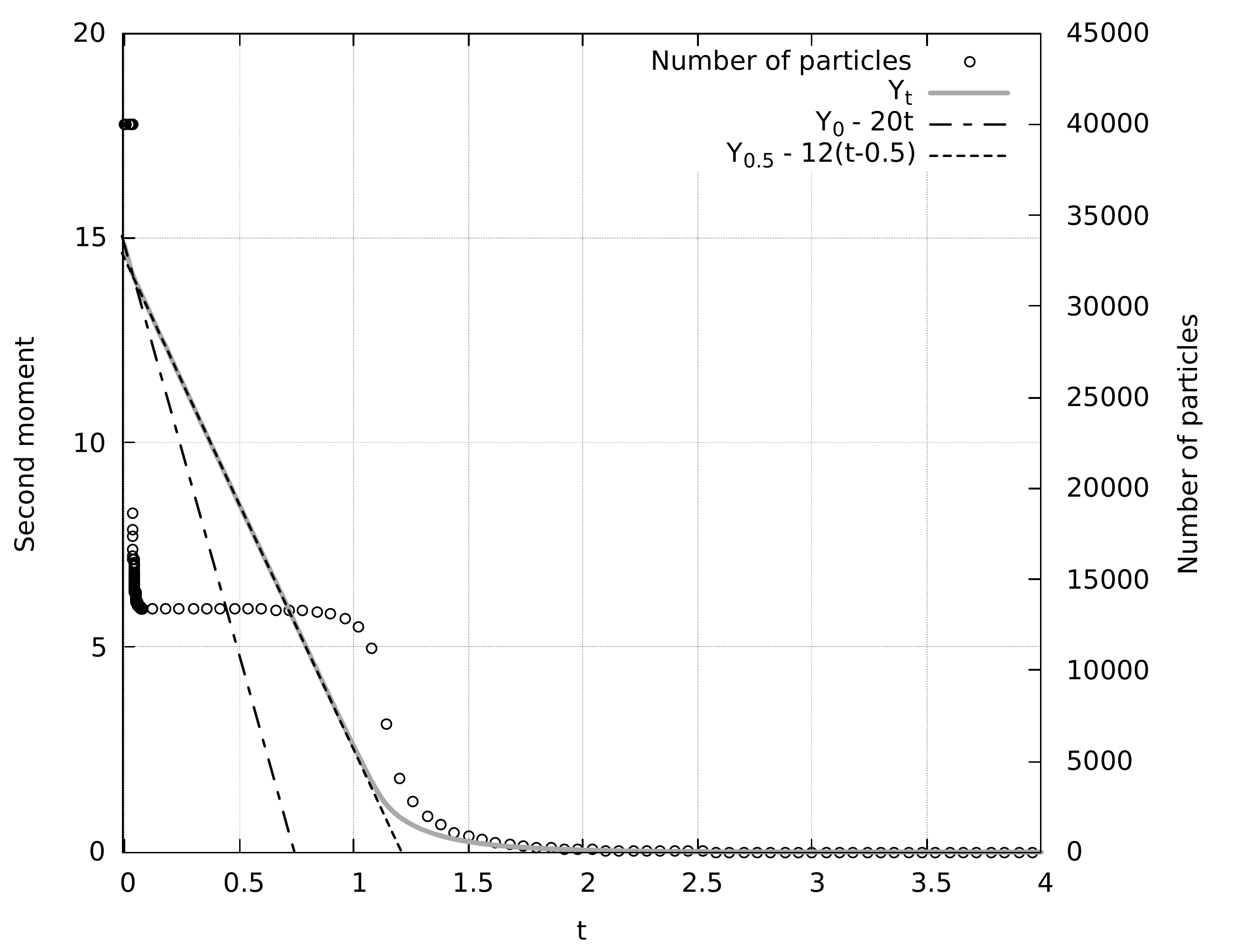}
	\end{center}
	\caption{We simulate $40 \times 10^3$ particles to approximate the system described in  \secref{sec:regularized_PKS_experiment} until it is fully singular. On this time scale, the continuous transfer of masses between the regular and singular component may be observed, by the curved second moment graph, and by the gradually decreasing graph of number of particles. The dashed lines correspond to the same ones as in \figref{fig:PKS_simulation_zoom}.}
	\label{fig:PKS_simulation_long_time}
\end{figure}

In particular, we consider the case $\chi = \mu = 1$, with total mass six times the critical mass, $M = 6\cdot 8 \pi$. We split the mass amongst a small bump function of mass $M_1 = 4 \cdot 8\pi$ supported on a disc of unit radius, which is separated far away from a bump function of mass $M_2 = 2 \cdot 8\pi$ that's supported on an ellipse with axes $1$ and $7$. These initial initial conditions are chosen so as to make the solution initially exhibit a linear decay of the second moment, then a sudden change of slopes due to the rapid formation of a singularity caused by the first bump function, and finally--a continuous deceleration, due to continual mass transfer from the lighter bump function to the formed atomic component. With the chosen parameters, the first two rates of change of the second moment should be $-20$ and $-12$. This can be observed in \figref{fig:PKS_simulation_zoom}. Further in time, the gradual transfer of mass may be seen as well, as shown in \figref{fig:PKS_simulation_long_time}.

\begin{figure}
	\begin{subfigure}{0.48\textwidth}
		\includegraphics[width=\linewidth]{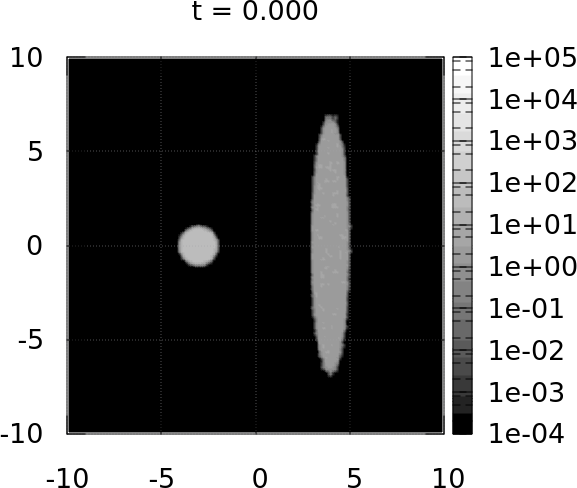}
	\end{subfigure}
	\begin{subfigure}{0.48\textwidth}
		\includegraphics[width=\linewidth]{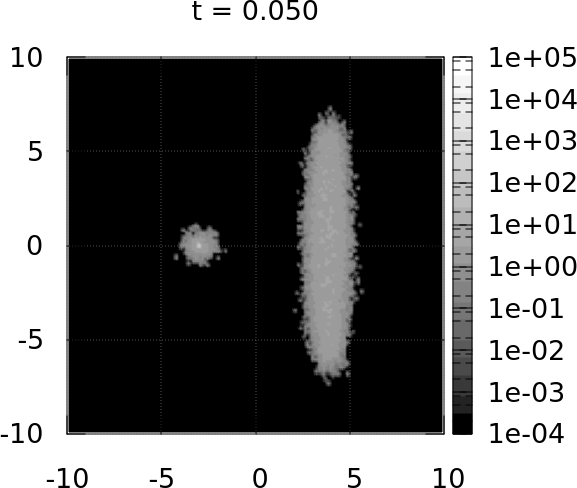}
	\end{subfigure}
	
	\vspace*{0.5cm}	
	
	\begin{subfigure}{0.48\textwidth}
		\includegraphics[width=\linewidth]{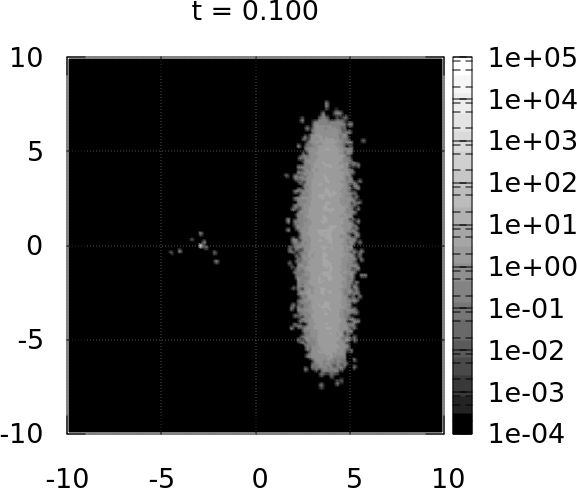}
	\end{subfigure}
	\begin{subfigure}{0.48\textwidth}
		\includegraphics[width=\linewidth]{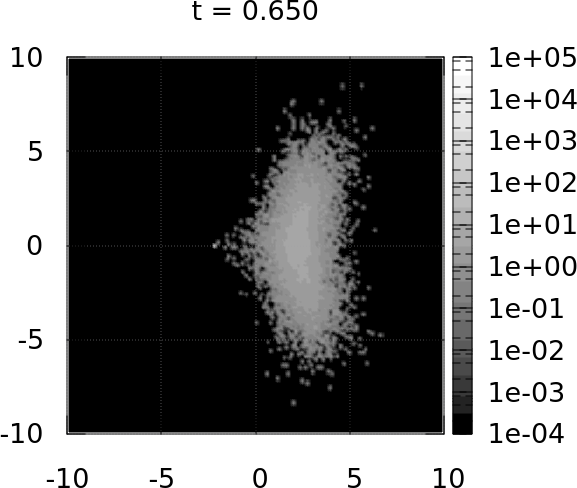}
	\end{subfigure}
	
	\vspace*{0.5cm}	
	
	\begin{subfigure}{0.48\textwidth}
		\includegraphics[width=\linewidth]{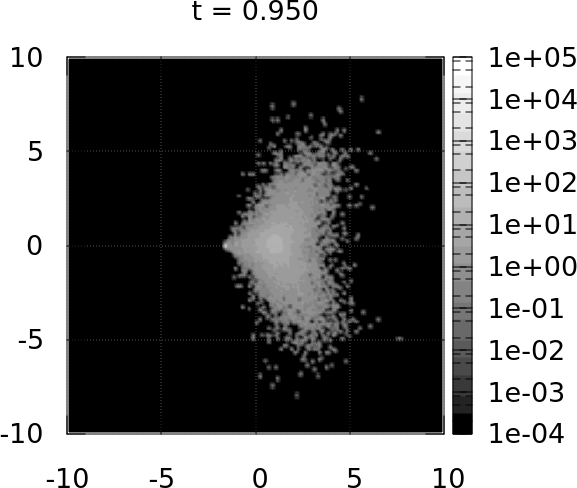}
	\end{subfigure}
	\begin{subfigure}{0.48\textwidth}
		\includegraphics[width=\linewidth]{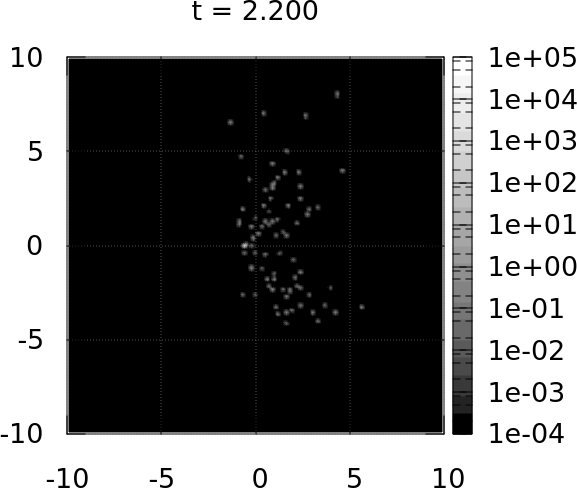}
	\end{subfigure}

	\caption{Snapshots of the interpolated mass density field $P_{ij}$ for the simulation described in \secref{sec:regularized_PKS_experiment}. The relation between this figure and \figref{fig:PKS_simulation_zoom} and \figref{fig:PKS_simulation_long_time} is given at the end of \secref{sec:mass_transfer_to_singularity}. All particles initially have the same mass.
	\label{fig:PKS_simulation_positions}}
\end{figure}

The underlying particle dynamics and collisions are illustrated in \figref{fig:PKS_simulation_positions}, where each snapshot corresponds to qualitatively different rates of change of the second moment in \figref{fig:PKS_simulation_zoom} and \figref{fig:PKS_simulation_long_time}: sudden mass coalescence of a tight aggregate (switch of slopes in \figref{fig:PKS_simulation_long_time}, and $t=0.050$ and $t=0.100$ in \figref{fig:PKS_simulation_positions}), attraction of mass without coalescence (linear decay in \figref{fig:PKS_simulation_zoom}, and $t=0.100$ and $t=0.650$ in \figref{fig:PKS_simulation_positions}), continuous slow and fast mass absorption (gradual deceleration in \figref{fig:PKS_simulation_long_time}, and $t=0.650$ and $t=0.950$ in \figref{fig:PKS_simulation_positions}), and the transformation of the PKS system to being essentially singular (flat part of the figure in \figref{fig:PKS_simulation_long_time}, and $t=2.200$ in \figref{fig:PKS_simulation_positions}).

\subsubsection{Interaction of singularities} 
\begin{figure}
\begin{center}
	\includegraphics[width=\linewidth]{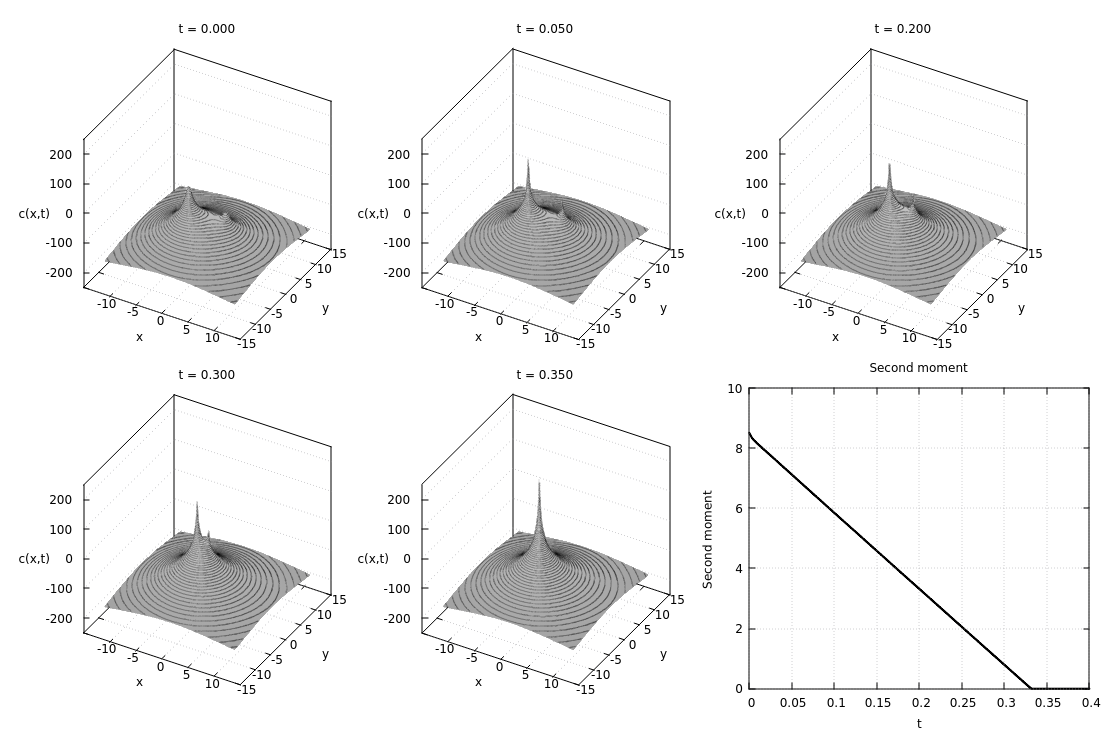}
\end{center}
	\caption{We set $\chi=\mu=1$ and initialize two small bumps functions at $(\pm 3, \pm 1)$ with supercritical masses $12\pi/5$ and $28\pi/5$. Each smooth bump quickly forms a singular component, and the $400 \times 10^3$ particle system reduces to a $\sim 2$ particle system. The formation and interaction of the singularities may be seen in the above snapshots of $c(x,t)$. After the initial formation of singularities, the second moment decreases linearly, as predicted by \eqref{eq:regularized_second_moment}. In this particular simulation, we used $[-15,15]^2$ as the computational domain, which we discretized using a $270 \times 270$ mesh, and set the time step to be $0.002$}
	\label{fig:interaction_of_singularities}
\end{figure}

In another experiment, we initialize a system in which \emph{two} singularities form and interact, as described in \figref{fig:interaction_of_singularities}. In this special case, the second moment is simply the square of the distance between the two singularities, the graph of which should be piecewise linear (as observed). We note that the numerical coalescence procedure avoids the ``washing out'' effect near the collision time in Figure 7 of \cite{fatkullin2013study}.

\subsection{Expanding MPKS with blow-up \label{sec:MPKS_experiment}} For the second example, we simulate blow-up with an increasing total second moment for the two species Keller-Segel system:
\begin{equation}
	\begin{cases}
		\partial_t \rho_1 &= \nabla \cdot (\mu_1 \nabla \rho_1 - \chi \rho_1 \nabla c ), \\
		\partial_t \rho_2 &= \nabla \cdot (\mu_2 \nabla \rho_2 - \chi \rho_2 \nabla c ), \\
		\Delta c &= -(\rho_1 + \rho_2),
	\end{cases}
\end{equation}
with $\int \rho_1 = M_1$ and $\int \rho_2 = M_2$. The interest in this phenomenon is described in \figref{fig:MPKS_blow_up_mass} and \secref{sec:MPKS_limit}. In particular, we show that when a two species PKS system is in this regime, the second moment of one component increases linearly, while the other decreases. Such semi-decoupled behavior was suggested in \cite{conca2011remarks}. We remark that the numerical method presented is well-suited for this investigation, as it can simulate the system in the entire plane.

We approximate this two  system using $N_0$ particles, the first $N_1 = \lfloor \eta_1 N_0 \rfloor$ of which have particle masses $M_1/N_1$, and distributed on the plane according to $\rho_1(\cdot, 0)$. Similarly, the last $N_2 = N_0 - N_1$ particles have masses $M_2/N_2$, and are distributed on the plane according to $\rho_2(\cdot, 0)$. Using \eqref{eq:nonuniform_masses_hydrodynamic} and \eqref{eq:nonuniform_masses_diffusion_coefficients}, we see that 
\begin{equation}
	\eta_i = \frac{M\mu_i}{M_i\mu},
\end{equation}
where
\begin{align}
	\label{eq:MPKS_mu}
	\mu = (M_1 + M_2)\left(\frac{\mu_1}{M_1} + \frac{\mu_2}{M_2}\right).
\end{align}
The particle system's diffusion coefficient $\tilde{\mu}$ is then
\begin{equation}
	\label{eq:2MPKS_diffusion_coefficient}
	\tilde{\mu} = \frac{\mu (M_1 + M_2)}{N_0}.
\end{equation}
Thus, for a two species MPKS system with component masses $M_1, M_2$ and diffusion coefficients $\mu_1, \mu_2$, we associate an $N_0$ particle system with two different possible particle masses. The diffusion coefficient for \eqref{eq:particle_dynamics} is given by \eqref{eq:2MPKS_diffusion_coefficient}. In this sense, the purpose of $\mu$ in \eqref{eq:MPKS_mu} is auxiliary.

When $\mu_1 > 2\mu_2$, it is always possible to choose component masses which will force a radially-symmetric system to blow-up with increasing second moment. In this case, $M_1^{max}$ and $M_2^{max}$ in \figref{fig:MPKS_blow_up_mass} can be shown to be
\begin{equation}
	\label{eq:MPKS_M_max}
	M_1^{max} = \frac{2\pi}{\chi} \cdot \frac{(\mu_1 - 2\mu_2)\mu_1}{\mu_1 - \mu_2}, \quad M_2^{max} = \frac{2\pi}{\chi} \cdot \frac{\mu_1^2}{\mu_1 - \mu_2}.
\end{equation}
For the experiments in this section, we simulate the two species system as described above, and choose the convenient parameters
\begin{equation}
\chi = 4,\quad \mu_1 = \frac{35}{2},\quad \mu_2 = \frac{35}{12}, \quad M_1 = 4, \quad M_2 = 24,
\end{equation}
which correspond to the auxiliary parameters
\begin{equation}
	\mu = 5, \quad \eta_1 = \frac{1}{2}, \quad \eta_2  = \frac{1}{2}.
\end{equation}

For the above masses, we consider three different initial conditions. Each respective solution exhibits linear growth in the first component's second moment, and decay in the second component's second moments, but at rates which depend on the initial distribution of mass. In particular, we choose the following initial conditions:

\begin{enumerate}
	\item \emph{Radially-symmetric component initial data.} We initialize both components as bump functions supported on a disc of radius $a = 0.35$ and centered at the origin.
	\item \emph{Non-symmetric component initial data.} We initialize the first component as a bump function of radius $a$ and centered at the origin, and the second component as a bump function supported on an ellipse centered at $(0.1, 0)$ with axes $2a$ and $a/2$, with the major axis parallel to the $y$-axis. 
	\item \emph{Component initial data on disjoint support.} We initialize each component on a bump function supported on a disc of radius $a$, where the first component is centered at $(a,-a)$, and the second at $(-a, a)$.
\end{enumerate}

\begin{figure}
\begin{center}
	\begin{subfigure}{0.48\textwidth}
		\includegraphics[width=\linewidth]{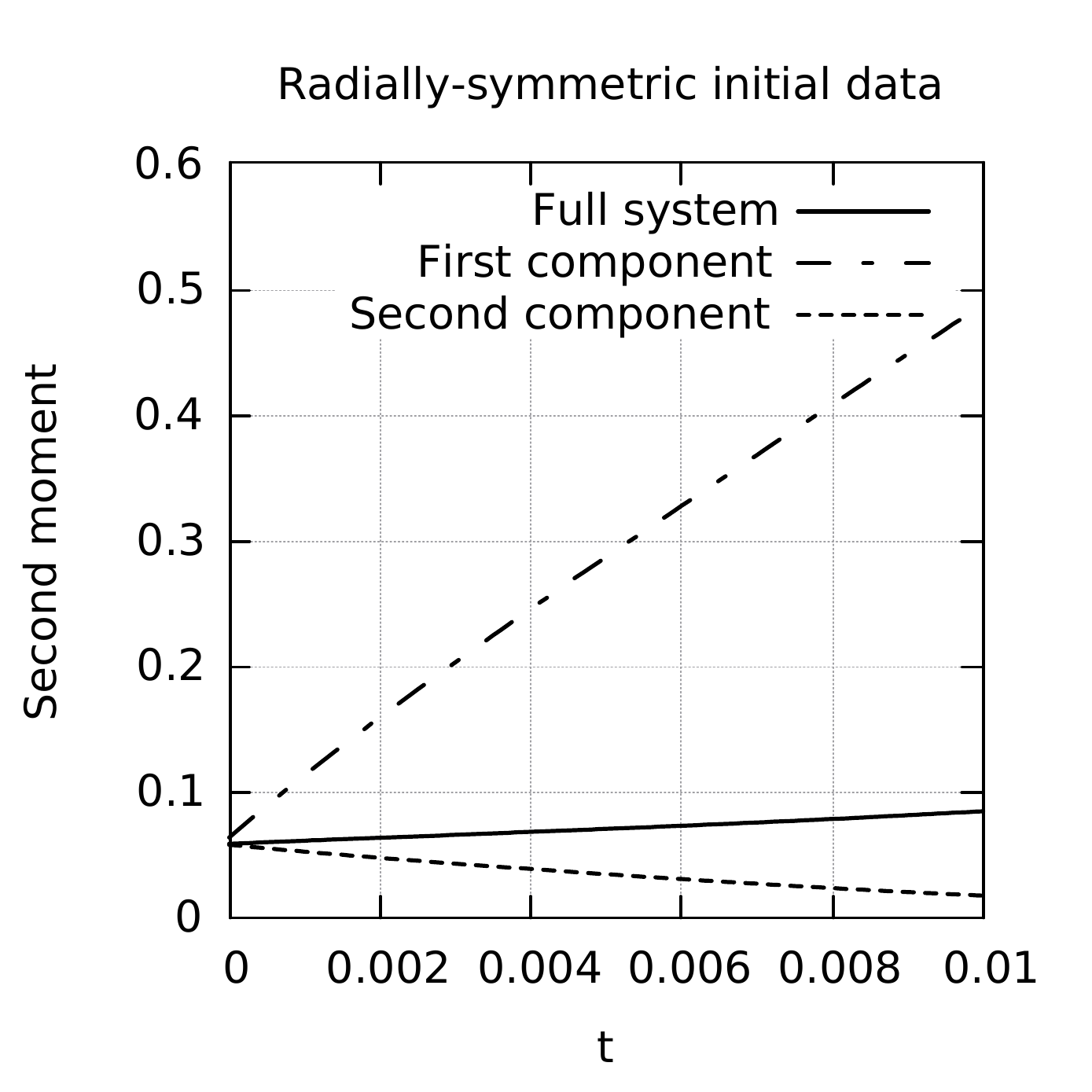}
	\end{subfigure}
	\begin{subfigure}{0.48\textwidth}
		\includegraphics[width=\linewidth]{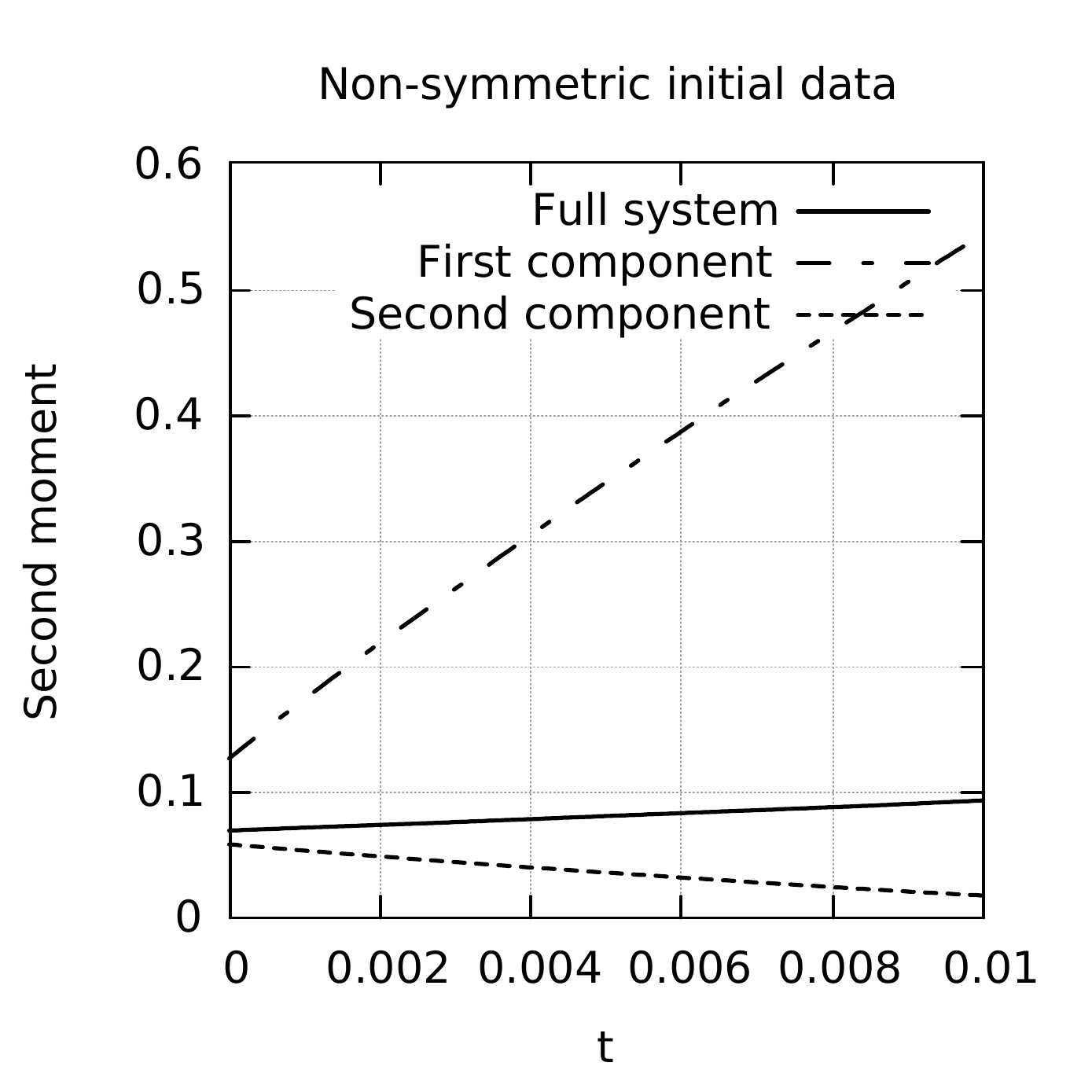}
	\end{subfigure}

	\begin{subfigure}{0.48\textwidth}
		\includegraphics[width=\linewidth]{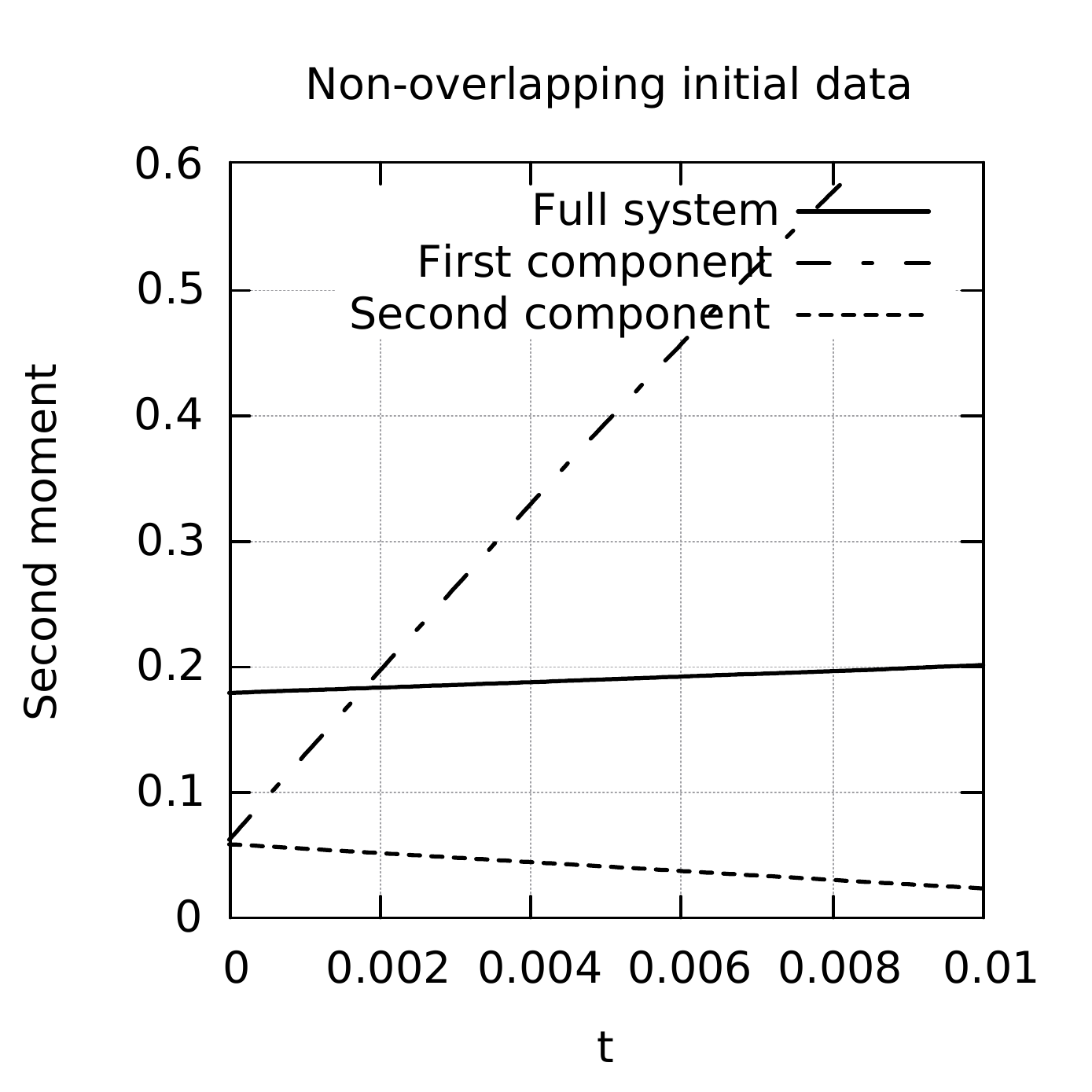}
	\end{subfigure}
	\begin{subfigure}{0.48\textwidth}
		\includegraphics[width=\linewidth]{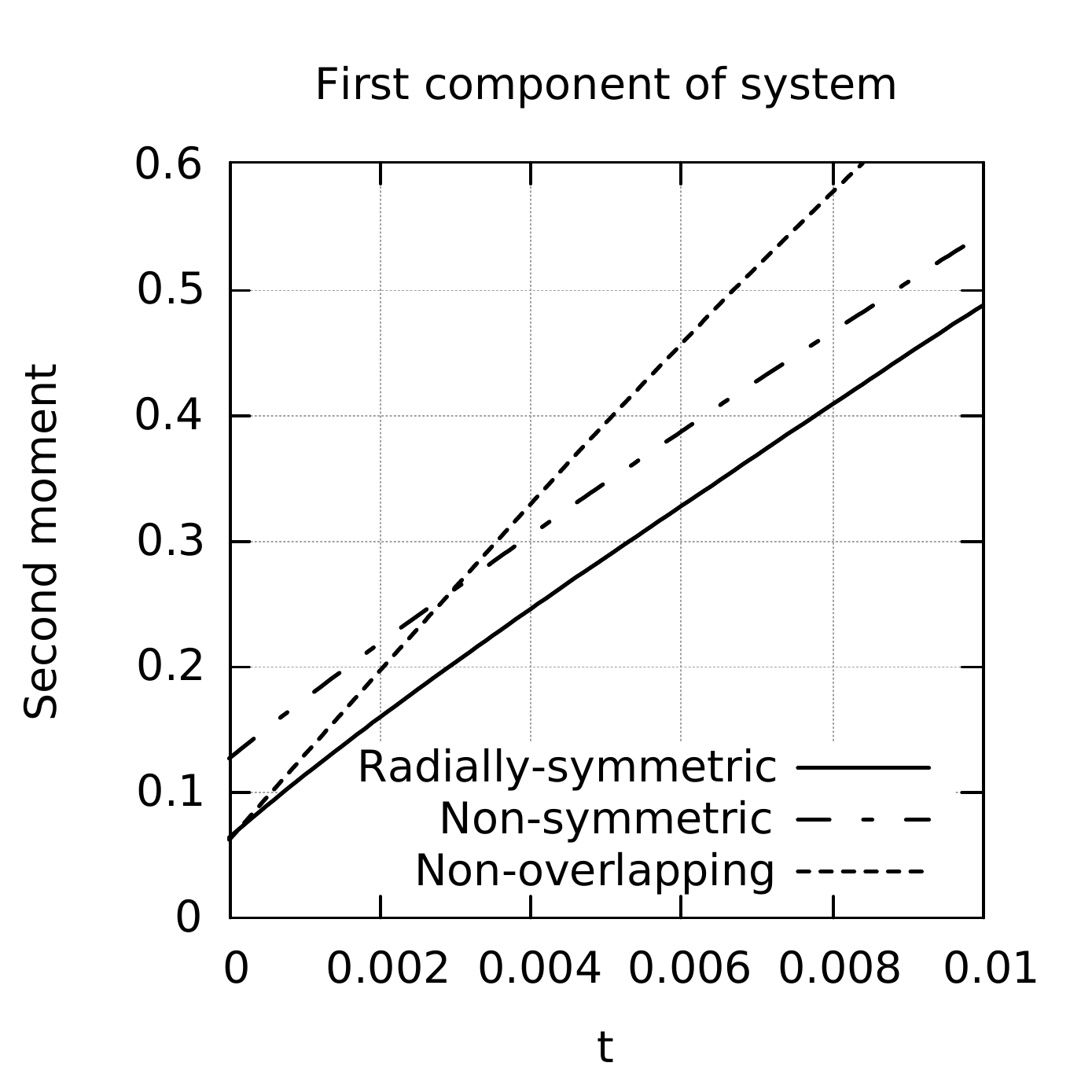}
	\end{subfigure}
\end{center}
	\caption{Evidence of the phenomenon described in \figref{fig:MPKS_blow_up_mass}, for initial conditions which are and are not radially-symmetric. As can be seen, the total second moment expands at a fixed rate, as do the individual components. However, the rate of change of the second moment of each component varies with the initial data. For these simulations, we used $10^6$ particles, and discretized $[-1.5, 1.5]$ using a $320 \times 320$ mesh for the computational grid. The initial conditions for each experiment are given in \ref{sec:MPKS_experiment}. We note that although the first component is expanding, there is evidence that it nonetheless blows up in the $L^\infty$ norm \cite{kurganov2014numerical}.}
	\label{fig:MPKS_second_moments}
\end{figure}

The results of these simulations can be seen in \figref{fig:MPKS_second_moments}. We note that although both components change linearly, their rates of change appear to depend on the initial conditions.

\section{Conclusion}
We investigated a planar particle system with nonuniform particle masses, in which particles interact via a logarithmically-singular kernel. As post-collision dynamics in such a system are undefined, we used the idea of particle coalescence in order to propagate the system further in time, and connected it to the theory of the squared Bessel process. We exploited this connection to develop an efficient numerical method for the simulation of the system, which has applications in the numerical approximation and regularization of a wide range of nonlinear Fokker-Planck equations, such as the multispecies Patlak-Keller-Segel model.

As mentioned before, properties of singularity formation in the MPKS are not fully understood, and have somewhat unexpected behavior, when compared to the PKS. For instance, singularities may form while the system's second moment is increasing. It would be interesting to further connect existing results with predicting a nonuniform particle system's behavior post-collision.

The question of coalescence in a system with memory arises naturally, as an analogue to the parabolic Keller-Segel model. In this case, the field $c(x,t)$ is replaced with the solution to the following equation,
\begin{equation}
	\partial_t c = \Delta c - k^2 c + \sum_{i} m_i \delta \left(x - X^{(i)}_t \right),
\end{equation}
which has the more biologically-meaningful intepretation of a chemoattractant which thermalizes at a finite rate, diffuses, decays, and is produced by the particles. This system will be investigated in future works.

\section{Acknowledgements}
This material is based upon work supported by the National Science Foundation under Grant DMS-1056471. The support is gratefully acknowledged.

The authors would like to thank Hailiang Liu for references related to the two species Keller-Segel model.

\bibliographystyle{siam}

\bibliography{PKS_bib}

\appendix
\section{Subtraction formula for indices\label{app:Subtraction_formula}}
If $\nu_i$ is the index of the full system described in \figref{fig:formed_cluster}, $\nu_f$ is the index of the same system after the particles inside the dashed lines coalesce, and $\overline{\nu}$ is the index of the subsystem inside the dashed line, then using \eqref{eq:Bessel_index} we have
\begin{align}
	\nu_i &= N - 2 - \frac{\chi}{8 \pi \tilde{\mu}} \left( M^2 - \sum_{j=1}^N m_j^2 \right), \\
	\nu_f &= N - N' - 1 - \frac{\chi}{8 \pi \tilde{\mu}} \left( M^2 - \sum_{j=N'+1}^N m_j^2 -(M')^2 \right), \\
	\overline{\nu} &= N' - 2 - \frac{\chi}{8 \pi \tilde{\mu}} \left( (M')^2 - \sum_{j=1}^{N'} m_j^2 \right),
\end{align}
from which it follows that
\begin{equation}
	\nu_f - \nu_i = -\left(\overline{\nu} + 1\right).
\end{equation}

\section{MPKS second moment\label{app:MPKS_second_moment_derivation}}
The evolution of the second moment of the MPKS can be computed as follows:
\begin{align}
	F'(t) &= \frac{d}{dt} \int_{\mathbb{R}^2} |x|^2 \sum_{i=1}^K \rho_i(x,t) dx \\
	&= \int_{\mathbb{R}^2} |x|^2 \sum_{i=1}^K \nabla \cdot ( \mu_i \nabla \rho_i - \chi \rho_i \nabla c) dx \\
	&= -2 \int_{\mathbb{R}^2} \sum_{i=1}^K \left( \mu_i \nabla \rho_i - \chi \rho_i \nabla c \right) \cdot x dx \\
	&= -2 \int_{\mathbb{R}^2} \sum_{i=1}^K \mu_i \nabla \rho_i \cdot x dx + 2\chi \int_{\mathbb{R}^2} \sum_{i=1}^K \rho_i \nabla c \cdot x dx \\
	&= 4 \int_{\mathbb{R}^2} \sum_{i=1}^K \mu_i \rho_i(x) dx \\
	\nonumber
	&\qquad \qquad - \frac{\chi}{\pi} \int_{\mathbb{R}^2 \times \mathbb{R}^2} \sum_{i,j=1}^K \rho_i(x) \rho_j(y) \frac{x-y}{|x-y|^2} dy \cdot x dx \\
	&= 4 \sum_{i=1}^K \mu_i M_i - \frac{\chi}{2\pi} \left( \int_{\mathbb{R}^2 \times \mathbb{R}^2} \sum_{i,j=1}^K \rho_i(x) \rho_j(y) \frac{x-y}{|x-y|^2} dy \cdot x dx \right. \\
	\nonumber	
	&\qquad \qquad  + \left. \int_{\mathbb{R}^2 \times \mathbb{R}^2} \sum_{i,j=1}^K \rho_i(y) \rho_j(x) \frac{y-x}{|y-x|^2} dx \cdot y dy \right) \\
	&= 4 \sum_{i=1}^K \mu_i M_i - \frac{\chi}{2\pi} \int_{\mathbb{R}^2} \sum_{i,j=1}^K \rho_i(y) \rho_j(x) dy dx \\
	&= 4 \sum_{i=1}^K \mu_i M_i - \frac{\chi M}{2\pi}\sum_{i=1}^K M_i \\
	&= \sum_{i=1}^K \left(4\mu_i - \frac{\chi M}{2\pi}\right)M_i.
\end{align}

\end{document}